\documentclass[11pt]{article} 
\pdfoutput=1

\usepackage{tikz}
\usetikzlibrary{positioning}
\usepackage[utf8]{inputenc} 
\usepackage{amsmath}
\usepackage{amssymb}
\usepackage{graphicx}
\usepackage{epsfig,epstopdf}
\usepackage{pdflscape}
\usepackage{rotating}

\usepackage{geometry} 
\usepackage{graphicx} 

\usepackage{booktabs} 
\usepackage{array} 
\usepackage{paralist} 
\usepackage{verbatim} 
\usepackage{subfig} 

\usepackage{fancyhdr} 
\usepackage{sectsty}
\allsectionsfont{\sffamily\mdseries\upshape} 

\usepackage[nottoc,notlof,notlot]{tocbibind} 
\usepackage[titles,subfigure]{tocloft} 

\title{Simplified three player Kuhn poker}
\author{J. Billingham\\ School of Mathematical Sciences, \\ The University of Nottingham, \\ Nottingham NG7 2RD, UK}

\begin{document}
\begin{abstract}
We study a very small three player poker game (one-third street Kuhn poker), and a simplified version of the game that is interesting because it has three distinct equilibrium solutions. For one-third street Kuhn poker, we are able to find all of the equilibrium solutions analytically. For large enough pot size, $P$, there is a degree of freedom in the solution that allows one player to transfer profit between the other two players without changing their own profit. This has potentially interesting consequences in repeated play of the game. We also show that in a simplified version of the game with $P>5$, there is one equilibrium solution if $5<P\leq P^* \equiv (5+\sqrt{73})/2$, and three distinct equilibrium solutions if $P>P^*$. This may be the simplest non-trivial multiplayer poker game with more than one distinct equilibrium solution and provides us with a test case for theories of dynamic strategy adjustment over multiple realisations of the game. 

We then study a third order system of ordinary differential equations that models the dynamics of three players who try to maximise their expectation by continuously varying their betting frequencies. We find that the dynamics of this system are oscillatory, with two distinct types of solution. Finally, we study a difference equation model, based on repeated play of the game, in which each player continually updates their estimates of the other players' betting frequencies. We find that the dynamics are noisy, but basically oscillatory for short enough estimation periods and slow enough frequency adjustments, but that the dynamics can be very different for other parameter values.
\end{abstract}
\maketitle

\section{Introduction}
Poker is a  multiplayer game of imperfect information. Although popular variants of poker, such as Texas Holdem and Omaha, are large, complex well-defined test problems for researchers in artificial intelligence \cite{Billings:2002:CP:512148.512158}, their size and complexity masks the fact that even small and simple toy poker games pose significant theoretical challenges. Most of the literature on poker, simplified or otherwise, focusses on solving the game, in the sense of finding Nash equilibrium solutions. Two player poker is a zero sum game, so all equilibrium solutions have the same expectation for each player, \cite{10.2307/j.ctt1r2gkx}. Recently, two player Limit Holdem, a large and complex game, was numerically solved in this sense \cite{15Science}, and a recently-developed, sophisticated AI that combines neural networks with on the fly equilibrium computation plays the even larger game of two player No Limit Holdem close to optimally \cite{DSAI}. However, even a simple two player, zero sum game leads to significant theoretical challenges when studied as a repeated game  (for example, two player Kuhn poker, also known as the AKQ game, \cite{Southey2009}).

Full street, three player Kuhn poker with pot size three units was introduced in \cite{Risk:2010:UCR:1838206.1838229} in order to test the performance of the counterfactual regret algorithm in three player games, a context in which convergence to an equilibrium solution is not guaranteed. Later, a family of equilibrium solutions was found analytically, \cite{Szafron:2013:PFE:2484920.2484962}, although it was not established whether other equilibrium solutions exist. This game was subsequently proposed as part of the 2015 Annual Computer Poker Competition to test the performance of AI players \cite{KuhnComp}. An indicator of the jump in dynamic complexity from two player to three player zero sum games is the fact that a repeated game for only three players with a deck of just four cards is a difficult, unsolved problem that attracts significant attention from the AI community. 

In this paper we study the one-third street version of three player Kuhn poker. Full and half street two player games are discussed extensively by Chen and Ankenman in \cite{MoP}. Following  their definitions, in the full street game Players 1, 2 and 3 can all be the first to bet, whilst in the one-third street game, Players 1 and 2 must check, and Player 3 makes the first decision (check or bet). It is this simpler game that we study here, but with the size of the initial pot, $P$ units, arbitrary. The natural choice, $P=3$, (each player contributes one unit to the pot before play starts) has been used in previous work, \cite{Risk:2010:UCR:1838206.1838229}, \cite{Szafron:2013:PFE:2484920.2484962}. In Section~\ref{sec_Kuhn} (details in Appendix~\ref{AKQJ_ap}), we determine all of the equilibrium solutions of three player, one-third street Kuhn poker. The equilibrium is unique and trivial for $P\leq 2$. For $P>2$, there is no unique equilibrium solution, but the nature of this non-uniqueness varies with $P$. When $P=5$ or $P\geq P^* \equiv \frac{1}{2}\left(5+\sqrt{73}\right) \approx 6.77$, there is a degree of freedom in the solution that allows one player to transfer profit between the other two players without changing their own profit.  This has potentially interesting consequences in repeated play of the game. For all other values of $P>2$, Player 3 bluffs with J and/or Q at equilibrium with a well-defined total bluffing frequency.

Our analysis of this game suggests that it is also of interest to study a simplified variant (Player 3 must check with K and Q, Player 2 must call a single bet with K, $P\geq 5$), which forms the subject of the rest of this paper. In Section~\ref{sec_simKuhn}, we show (details in Appendix~\ref{restAKQJ_ap}) that this simplified game has one equilibrium solution when $5<P\leq P^*$, and three distinct equilibrium solutions when $P>P^*$. Moreover, for $P>7$ it is not clear which, if any, of these equilibria would be prefered in a game between rational players. One of our main aims in this paper is to introduce this simplified three player variant as perhaps the most straightforward nontrivial multiplayer poker game with more than one distinct equilibrium solution. It contains just three strategic decisions, characterized by three betting frequencies, one for each player. This suggests several ways in which the repeated game can be modelled. In section~\ref{sec_ode}, we study a third order system of ordinary differential equations that models the dynamics of three players who try to maximise their expectation by continuously varying their betting frequencies. We find that the dynamics of this sytem are oscillatory, with two distinct types of solution. In section~\ref{sec_diff}, we study a difference equation model, based on repeated play of the game, in which each player continually updates their estimates of the other players' betting frequencies. We find that the dynamics are noisy, but basically oscillatory for short enough estimation periods and slow enough frequency adjustments, but that the dynamics can be very different for other parameter values.

\section{One-third street, three player, Kuhn poker}\label{sec_Kuhn}
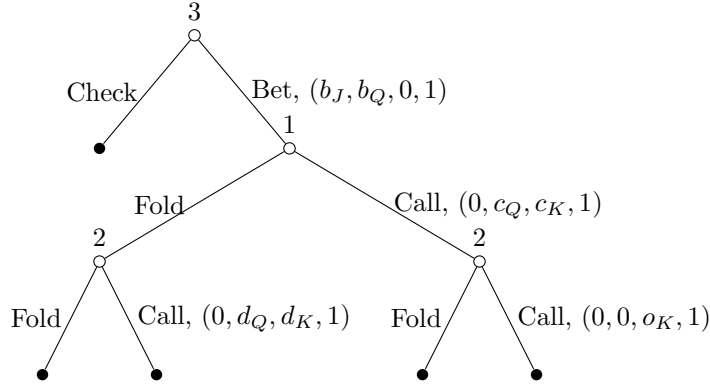
\begin{figure}
 	\begin{center}
    \small
    \begin{tikzpicture}[thin,
      level 1/.style={sibling distance=25mm},
      level 2/.style={sibling distance=50mm},
      level 3/.style={sibling distance=15mm},
      every circle node/.style={minimum size=1.5mm,inner sep=0mm}]

      \node[circle,draw,label=above:$3$] (root) {}
        child { node [circle,fill] {}
          edge from parent
            node[left] {Check}}
        child { node [circle,draw,label=above:$1$] {}
          child {
            node[circle,draw,label=above:$2$] (node-A) {}
              child {node [circle,fill] {}
                edge from parent
                  node[left] {Fold}}
              child {node [circle,fill] {}
                 edge from parent
                   node[right] {Call, $(0, d_Q, d_K, 1)$}}
              edge from parent
                node[left] {Fold}}
          child {
            node[circle,draw,label=above:$2$] (node-B) {}
              child {node [circle,fill] {}
                edge from parent
                  node[left] {Fold}}
              child {node [circle,fill] {}
                 edge from parent
                   node[right] {Call, $(0, 0, o_K, 1)$}}
              edge from parent
                node[right] {Call, $(0, c_Q, c_K, 1)$}}
           edge from parent
             node[right] {Bet, $(b_J, b_Q, 0, 1)$}};
    \end{tikzpicture}
    \end{center}
    \caption{The decision tree for three player, one-third street Kuhn poker. Open circles are decision nodes (labelled by the player making the decision), whilst solid circles are terminal nodes. Betting and calling frequencies with $(J, Q, K, A)$ are also shown.\label{fig_tree}}
  \end{figure}

The deck contains four cards, A $>$ K $>$ Q $>$ J, and each player is dealt a single card at random. There is a pot of $P$ units. Players 1 and 2 are forced to check. Player 3 can then either check, in which case there is a showdown, or bet one unit. If Player 3 bets, Player 1 must either call or fold. Player 2 must then either call (overcall if Player 1 has called) or fold. If Player 1 and/or Player 2 calls, there is a showdown at which the player with the best card wins the pot and all the bets, otherwise Player 3 wins the pot. This is the one-third street version of three player Kuhn poker, which was introduced in \cite{Risk:2010:UCR:1838206.1838229}. Figure~\ref{fig_tree} shows the decision tree\footnote{Note that this is not the full game tree as it does not show the hidden cards (information sets) or the payoffs at the terminal nodes.}.

Although each player is free to make either of their available decisions with any card, assuming rational players:
\begin{enumerate}
\item  Player 3 will bet with A, and Players 1 and 2 will call or overcall with A (A always wins at showdown).
    \item Players 1 and 2 will always fold J when Player 3 bets (J always loses at showdown).
        \item Player 2 will not overcall with Q after a call from Player 1 (Q cannot be the best card).
    \item Player 3 will not bet with K (Either Player 1 or Player 2 will hold A two-thirds of the time, call and win, which leads to a loss that outweighs any profit Player 3 may make from Player 1 or 2 calling with Q).
\end{enumerate}
As an aside, note that point 4. above assumes a level of rationality that exceeds that of many recreational poker players. Betting in situations in which no worse hand will call and no better hand will fold is endemic amongst weak human players, \cite{EMPP}. 

The strategy parameters that we need to consider are:
  \begin{itemize}
  \item[{\em Player 3}:] Bluffing frequencies $b_J$ and $b_Q$ with J and Q.
  \item[{\em Player 1}:] Calling frequencies $c_Q$ and $c_K$ with Q and K.
  \item[{\em Player 2a}:] Calling frequencies $d_Q$ and $d_K$ with Q and K after Player 2 folds.
  \item[{\em Player 2b}:] Overcalling frequency $o_K$ with K after Player 2 calls.
  \end{itemize}

When $P < 2$, $b_J = b_Q = c_Q = c_K = d_Q = d_K=0$, with $o_K$ undetermined, is the unique, trivial equilibrium solution. Unless the pot is large enough, Player 3 has no incentive to bluff, and Players 1 and 2 therefore have no incentive to call. We will now assume that $P\geq2$. We show in Appendix~\ref{AKQJ_ap}, and confirm by symbolic algebra in Appendix~\ref{sec_Math2}, that the equilibrium solutions are
 \begin{equation}
 0\leq b_J+b_Q \leq \frac{2}{3},~~c_Q = c_K = d_Q = d_K =o_K=0~~\mbox{for $P=2$,}\label{eq1}
 \end{equation}
 \begin{equation}
 b_J+b_Q = \frac{2}{P+1},~~c_Q = c_K = 0,~~d_Q = 0,~~d_K = \frac{2P-4}{P+1},~~o_K=0\label{eq2}
 \end{equation}
\[
\mbox{for $2 < P < 5$,}
\]
 \begin{equation}
 \frac{1}{3}\leq b_J+b_Q \leq \frac{2}{5},~~c_Q = c_K = 0,~~d_Q = 0,~~d_K = 1,~~o_K=0\label{eq3}
 \end{equation}
\[\mbox{for $P= 5$,}
\]
 \begin{equation}
 b_J+b_Q = \frac{2}{P},~~b_Q>\frac{12+5P-P^2}{6P(P+1)},~~c_Q = 0,~~c_K = \frac{P-5}{P+1},~~d_Q = 0,~~d_K = 1,~~o_K=0\label{eq4a}
 \end{equation}
\[\mbox{for $5<P<P^*\equiv \frac{1}{2}(5+\sqrt{73}) \approx 6.67$,}
\]
 \begin{equation}
 b_J+b_Q = \frac{2}{P},~~c_Q = 0,~~c_K = \frac{P-5}{P+1},~~d_Q = 0,~~d_K = 1,~~o_K=0\label{eq4}
 \end{equation}
\[\mbox{for $P\geq P^*$,~~(Solution A)}
\]
 \begin{equation}
 b_J = \frac{2}{P},~~b_Q=0,~~0\leq c_Q \leq \frac{2}{P+4},~~c_K = \frac{P-5}{P+1},~~d_Q = 0,~~d_K = 1,~~o_K=0\label{eq5}
 \end{equation}
\[\mbox{for $P\geq P^*$.~~(Solution B)}
\]
Note that Player 2 neither calls with Q nor overcalls with K at equilibrium ($d_Q = 0$ and $o_K = 0$), which is not obvious {\it a priori}. Note also that $P^* \equiv \frac{1}{2}\left(5+\sqrt{73}\right)$ is the only positive root of $12+5P-P^2 = 0$, so the constraint on $b_Q$ in (\ref{eq4a}) is satisfied by all positive $b_Q$ when $P\geq P^*$.

\begin{figure}
 \begin{center}
 \includegraphics[width=\textwidth]{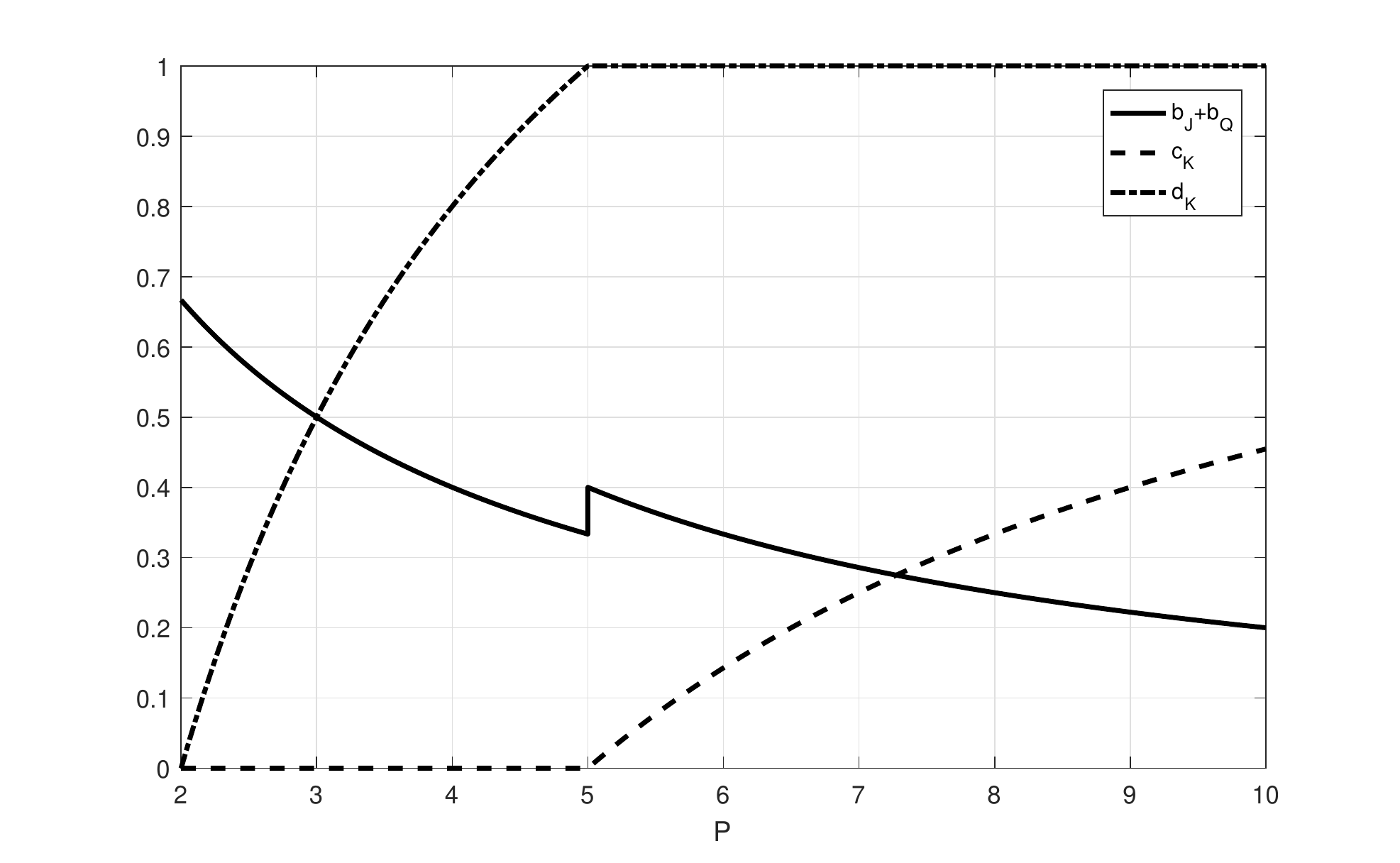}
 \caption{The equilibrium bluffing and calling frequencies for one-third street, three player Kuhn poker.}\label{fig_f}
 \end{center}
 \end{figure}

At the bifurcation points $P=2$ and $P=5$, a range of total bluffing frequencies, $b_J+b_Q$ is possible for Player 3. Away from these two points, there is a fixed bluffing frequency for Player 3. However, with the exception of Solution B, the choice of J or Q as a bluffing card for Player 3 is not constrained; only $b_J+b_Q$ is prescribed. This indeterminacy arises because Players 1 and 2 never call with Q at equilibrium. In Solution B, which only exists for $P>P^*\equiv \frac{1}{2}(5+\sqrt{73}) \approx 6.77$, Player 1 can call with Q with a frequency up to $2/(P+4)$ without affecting her expectation, and Players 2 or 3 cannot exploit her, neither by overcalling with K (Player 2) nor by bluffing more with J (Player 3), provided that Player 3 bluffs only with J ($b_J = 2/P$, $b_Q = 0$). The equilibrium frequencies are plotted in Figure~\ref{fig_f}. For $2 < P < 5$, Player 2 uses K to catch Player 3's bluffs. For $P>5$, Player 2 cannot call often enough with K, and Player 1 must call with K at frequency $(P-5)/(P+1)$, which allows Player 3 to bluff slightly more often.

\begin{figure}
 \begin{center}
 \includegraphics[width=\textwidth]{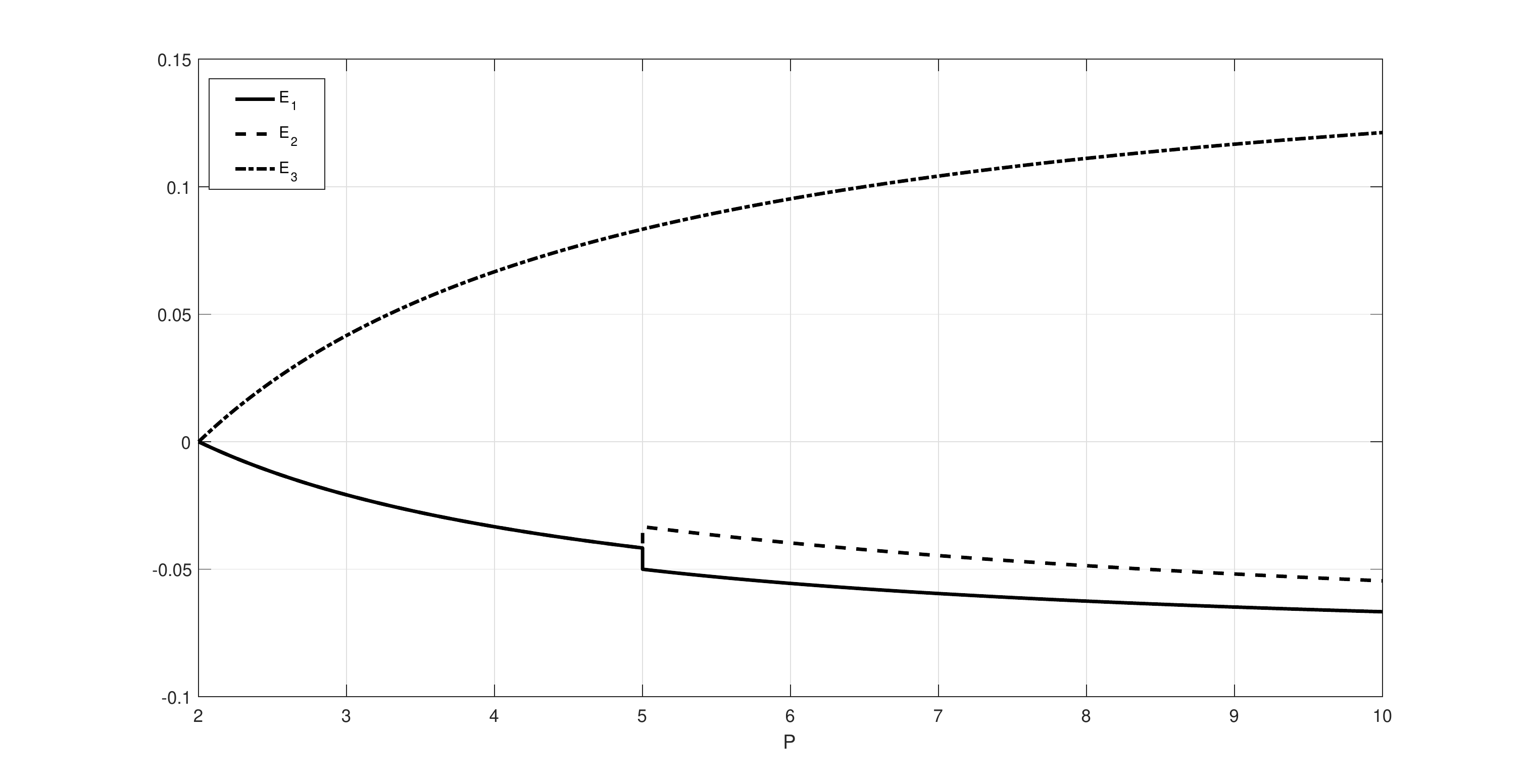}
 \caption{The ex-showdown expectations for one-third street, three player Kuhn poker with $c_Q=0$ (i.e. not solution B).}\label{fig_E}
 \end{center}
 \end{figure}

 The ex-showdown expectations of these solutions are
 \begin{equation}
 E_1 = E_2 = E_3 = 0~~\mbox{for $P\leq2$,}
 \end{equation}
 \begin{equation}
 E_1 = -\frac{P-2}{12(P+1)},~~E_2 = -\frac{P-2}{12(P+1)},~~E_3 = \frac{P-2}{6(P+1)}~~\mbox{for $2 \leq P < 5$,}
 \end{equation}
 \begin{equation}
 E_1 = -\frac{1}{8}(b_J+b_Q),~~E_2 = -\frac{1}{12} + \frac{1}{8}(b_J+b_Q),~~E_3 = \frac{1}{12}~~\mbox{for $P=5$,}
 \end{equation}
 \begin{equation}
 E_1 = -\frac{P-2}{12P},~~E_2 = -\frac{(P-1)(P-2)}{12P(P+1)},~~E_3 = \frac{P-2}{6(P+1)}~~\mbox{for $P>5$, (Solution A)}
 \end{equation}
\begin{equation}
E_1 = -\frac{P-2}{12P},~~E_2 = -c_Q -\frac{(P-1)(P-2)}{12P(P+1)},~~E_3 = c_Q+\frac{P-2}{6(P+1)}~~\mbox{for $P \geq P^*$, (Solution B)}
\end{equation}
 and are plotted in Figure~\ref{fig_E} with $c_Q = 0$. For all $P$, Player 3 has the chance to check with K and thereby realize its potential at showdown (poker players say that Player 3 has position). In contrast, at equilibrium Players 1 and 2 sometimes fold K when it is the best card and would win the pot at showdown. This means that Player 3's ex-showdown expectation, $E_3$, is positive, whilst those of Players 1 and 2, $E_1$ and $E_2$, are negative.

Note that Player 3's expectation, $E_3$, is a continuous function of $P$ at equilibrium, whilst those of the other players, $E_1$ and $E_2$, are discontinuous at $P=5$, with Player 1 losing more than Player 2 for $P>5$. At $P=5$, Player 3 can transfer profit between Players 1 and 2 at equilibrium by varying her total bluffing frequency  ($\frac{1}{3} \leq b_J+b_Q \leq \frac{2}{5}$). Similarly, when $P\geq P^*$, Player 1 can transfer profit from Player 2 to Player 3 by choosing $0\leq c_Q \leq \frac{2}{P+4}$. Player 3 should therefore choose $b_J = 2/P$ for $P \geq P^*$ in order to maximise her potential profit. We now discuss some possible implications of this for repeated play of the game.

\subsection{Repeated one-third street Kuhn poker}
If $P=5$, Player 3 can transfer profit between Players 1 and 2 at equilibrium by appropriately choosing her bluffing frequency, $b_J+b_Q$. If $P>P^*$, Player 1 can transfer profit between Players 2 and 3 by changing her calling frequency with Q, $c_Q$. In either case, in repeated play of the game with rotation, if two of the players form an alliance they can transfer profit to each other. 

As an example, consider the game with $P=5$. Player 3 controls the distribution of profit between Players 1 and 2 at equilibrium. If players 3 and 1 are in an alliance, Player 1 will use $c_K = c_Q = 0$, and if Player 3 chooses $b_J = 1/3$, $b_Q = 0$, the expectations of the three players are independent of Player 2's choice of strategy and Player 1's expectation is maximised. However, if Players 3 and 2 are in an alliance (perhaps the same two players after rotation of seats), with Player 2 choosing $d_Q = 0$, $d_K=1$ and $o_K = 0$ and Player 3 choosing $b_J+b_Q = 2/5$, the expectations are
\[
E_1 = \frac{5}{24} c_Q \left(b_J-\frac{2}{5}\right) - \frac{1}{20},~~E_2 = - \frac{5}{24} c_Q \left(b_J-\frac{1}{5}\right) - \frac{1}{30} + \frac{1}{60} c_K,~~E_3 = \frac{1}{12} - \frac{1}{60}c_K +\frac{1}{24} c_Q.
\]
By choosing non-zero values of $c_Q$ and/or $c_K$ (i.e. calling some fraction of the time with Q and/or K), Player 1 can transfer profit between Players 2 and 3. Note that by choosing $b_J < 1/5$, Player 3 can protect Player 2 from this possibility, but that if she chooses $b_J>1/5$, Player 1 can choose to target either Player 2 by calling with Q (at some cost to herself) or Player 3 by calling with K (at no cost to herself). Furthermore, as the players' positions rotate, the player who is not in the alliance will be Player 3 one-third of the time, which gives her an opportunity to decide how to treat the other two players. 

A further complicating factor is that the strategy parameters are betting and calling frequencies. The actual game is played with the cards dealt, not with publicly available values of the parameters, which must be estimated by the players from the information that they have. Game theory usually assumes rationality in the players, but it is not clear whether the level of rationality assumed when asserting that a player will not fold with A is the same as that assumed when asserting that a player will estimate another player's frequency of calling with Q when they are Player 1 and exploit any opportunity for profit that this might present. This is where skill enters the picture and an AI player should  be able to achieve win rates that exceed those of a rational but imperfect human player, or indeed an inferior AI. 

In \cite{Szafron:2013:PFE:2484920.2484962}, a similar possibility is discussed for an equilibrium solution of the full street version of the game. Note that there is numerical evidence from CFR solutions of the full street game with a range of pot sizes, \cite{JBCFR}, that $P=3$ (the pot size used in \cite{Szafron:2013:PFE:2484920.2484962}) is a bifurcation point (along with $P=4$ and $P=5$) in the same way that $P=5$ is a bifurcation point in the one-third street game studied here. It is not clear how much of the complicated equilibrium solution studied in \cite{Szafron:2013:PFE:2484920.2484962} remains when $ P \ne 3$, and whether there is any analogue of the continuous range of solutions of this type that exist in the one-third street game for $P > P^*$. 

\section{Simplified one-third street, three player, Kuhn poker (SKP)}\label{sec_simKuhn}
In this section we consider simplified one-third street three player Kuhn poker (SKP), in which:
\begin{itemize}
\item $P > 5$,
\item Player 3 must check with Q ($b_Q \equiv 0 $),
\item Player 1 cannot call with Q ($c_Q \equiv 0$),
\item Player 2 must call with K ($d_K \equiv 1$).
\end{itemize}
An immediate consequence of these restrictions is that Player 2 never overcalls with K when Player 1 calls (since Player 1 must then have A). Note that in the full game, which we analysed in the previous section, when $P>5$ the equilibrium solution has $d_K = 1$ and $o_K = 0$, and also allows $c_Q = 0$. The main simplification we make, which leads to the existence of multiple, distinct equilibrium solutions, is that Player 3 is not allowed to bluff with Q ($b_Q=0$).

\begin{figure}
 	\begin{center}
    \small
    \begin{tikzpicture}[thin,
      level 1/.style={sibling distance=25mm},
      level 2/.style={sibling distance=50mm},
      level 3/.style={sibling distance=15mm},
      every circle node/.style={minimum size=1.5mm,inner sep=0mm}]

      \node[circle,draw,label=above:$3$] (root) {}
        child { node [circle,fill] {}
          edge from parent
            node[left] {Check}}
        child { node [circle,draw,label=above:$1$] {}
          child {
            node[circle,draw,label=above:$2$] (node-A) {}
              child {node [circle,fill] {}
                edge from parent
                  node[left] {Fold}}
              child {node [circle,fill] {}
                 edge from parent
                   node[right] {Call, $(0, d_Q, 1, 1)$}}
              edge from parent
                node[left] {Fold}}
          child {
            node[circle,draw,label=above:$2$] (node-B) {}
              child {node [circle,fill] {}
                edge from parent
                  node[left] {Fold}}
              child {node [circle,fill] {}
                 edge from parent
                   node[right] {Call, $(0, 0, 0, 1)$}}
              edge from parent
                node[right] {Call, $(0, 0, c_K, 1)$}}
           edge from parent
             node[right] {Bet, $(b_J, 0, 0, 1)$}};
    \end{tikzpicture}
    \end{center}
    \caption{The decision tree for simplified three player, one-third street Kuhn poker. Open circles are decision nodes (labelled by the player making the decision), whilst solid circles are terminal nodes. Betting and calling frequencies with $(J, Q, K, A)$ are also shown.\label{fig_tree2}}
  \end{figure}
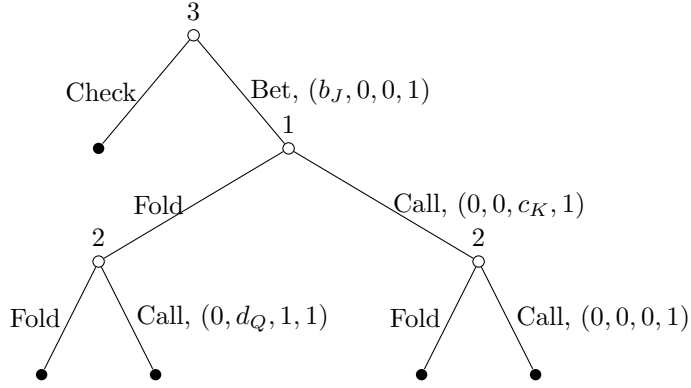

The remaining nontrivial decisions are:
\begin{enumerate}
\item[{\em Player 3}:] bluffing frequency with J, $b_J$,
\item[{\em Player 1}:] calling frequency with K, $c_K$,
\item[{\em Player 2}:] calling frequency with Q, $d_Q$.
\end{enumerate}
The simplified betting tree is shown in Figure~\ref{fig_tree2}.

\begin{figure}
 \begin{center}
 \includegraphics[width=\textwidth]{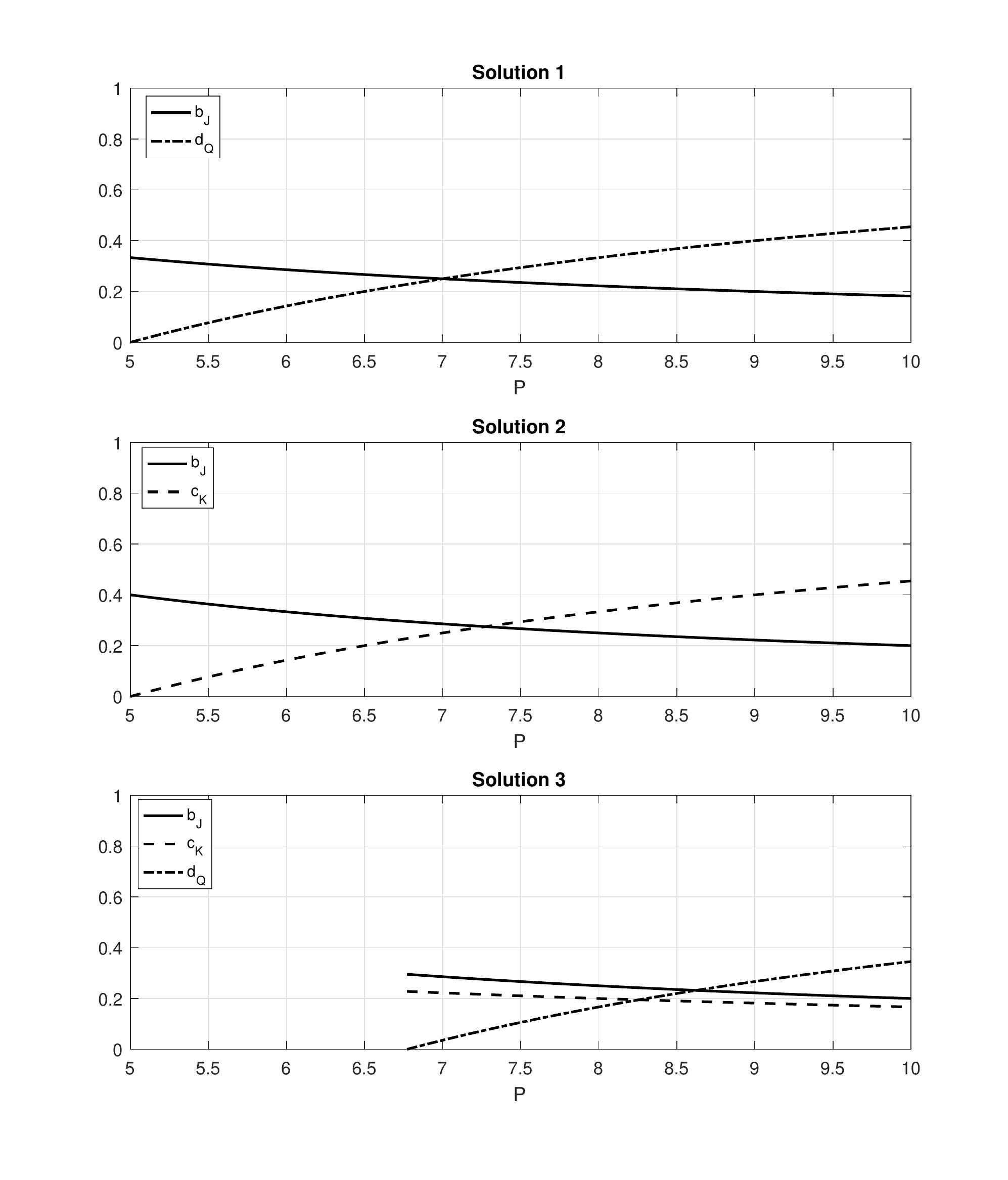}
 \caption{The three equilibrium solutions of simplified, one-third street, three player Kuhn poker.}\label{fig_resf}
 \end{center}
 \end{figure}

  We demonstrate  in Appendix~\ref{restAKQJ_ap}, and confirm numerically in Appendix~\ref{sec_Math1}, that the equilibrium solutions are
  \begin{equation}
  b_J = \frac{2}{P+1},~~c_K = 0,~~d_Q = \frac{P-5}{P+1},~~\mbox{(Solution 1)}\label{eqn_rsol1}
  \end{equation}
  \begin{equation}
  b_J = \frac{2}{P},~~c_K = \frac{P-5}{P+1},~~d_Q = 0,~~\mbox{(Solution 2)}\label{eqn_rsol2}
  \end{equation}
  \[~~\mbox{for $P \geq P^* = \frac{1}{2}(5+\sqrt{73}) \approx 6.77$,}\]
  \begin{equation}
  b_J = \frac{2}{P},~~c_K = \frac{2}{P+2},~~d_Q = \frac{P^2-5P-12}{P(P+1)},~~\mbox{(Solution 3)}\label{eqn_rsol3}
  \end{equation}
  \[~~\mbox{for $P \geq P^*$.}\]
In Solution 1, Player 2 calls with Q at a non-zero frequency, whilst Player 1 folds K. In Solution 2, Player 1 calls with K at the same non-zero frequency, whilst Player 2 folds Q (this is also a solution in the unsimplified, one-third street game). In Solution 3, Players 1 and 2 call with K and Q respectively at a non-zero frequency. These frequencies are illustrated in Figure~\ref{fig_resf}.

\begin{figure}
 \begin{center}
 \includegraphics[width=\textwidth]{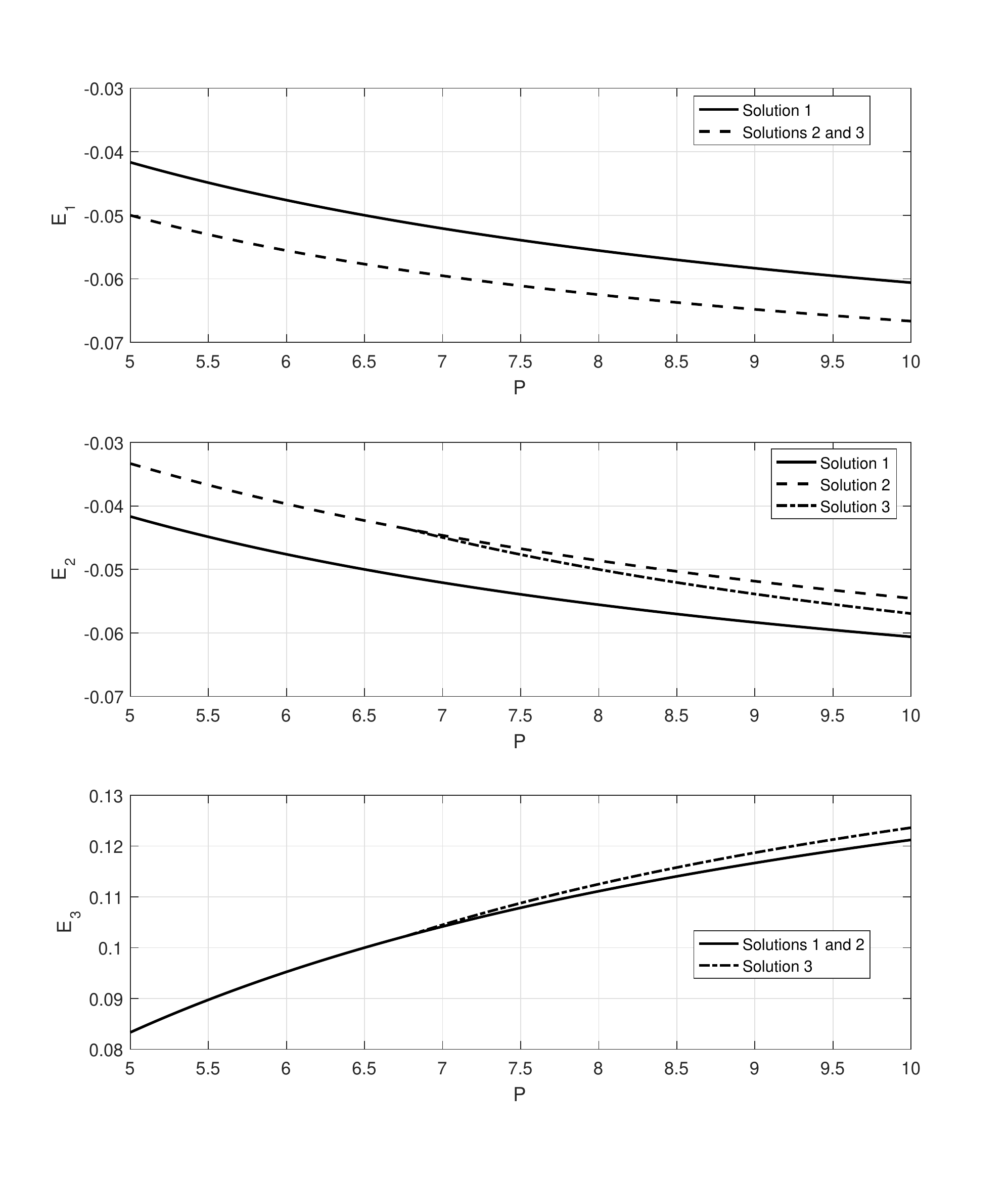} 
 \caption{The expectations of the three equilibrium solutions of simplified, one-third street, three player Kuhn poker.}\label{fig_resE}
 \end{center}
 \end{figure}

 The ex-showdown expectations,
 \begin{equation}
 E_1 = -\frac{P-2}{12(P+1)},~~E_2 = -\frac{P-2}{12(P+1)},~~E_3 = \frac{P-2}{6(P+1)},~~\mbox{(Solution 1)}\label{eqn_re1}
 \end{equation}
 \begin{equation}
 E_1 = -\frac{P-2}{12P},~~E_2 = -\frac{(P-1)(P-2)}{12P(P+1)},~~E_3 = \frac{P-2}{6(P+1)},~~\mbox{(Solution 2)}\label{eqn_re2}
 \end{equation}
 \begin{equation}
 E_1 = -\frac{P-2}{12P},~~E_2 = -\frac{P^2-P-8}{12P(P+2)},~~E_3 = \frac{2P^2-P-12}{12P(P+2)},~~\mbox{(Solution 3)}\label{eqn_re3}
 \end{equation}
are illustrated in Figure~\ref{fig_resE}.
Note that  Player 1 has her greatest expectation in Solution 1 and Player 2 in Solution 2. Player 3 has her greatest expectation in Solution 3 if $P > P^*$, but has the same expectation in each of Solutions 1 and 2 if $P\leq P^*$. However, if each player is restricted to choosing one of the three possible options (effectively two options for Player 3), we find that the maxmin strategy (the best worst-case payoff) is:
\begin{itemize}
\item Player 1: Solution 1,
\item Player 2: Solution 1 for $P<7$, Solution 2 for $P>7$,
\item Player 3: Solution 1.
\end{itemize}
This was determined by numerical computation. We conclude that for $P<7$, Solution 1 is a rational choice for the three players. However, for $P>7$, this suggests that Player 1 has an incentive to choose Solution 1 ($c_K = 0$) and Player 2 to choose Solution 2 ($d_Q = 0$). Although the maxmin strategy for Player 3 under these constraints is Solution 1, $c_K=d_Q=0$ is not part of an equilibrium solution and leaves Players 1 and 2 open to exploitation by Player 3 by increasing $b_J$. This strongly suggests that the outcome of a repeated game between rational players is unlikely to settle at an equilibrium solution if $P>7$, and, depending on the dynamics of the players' decision making, this may also be the case for $P\leq 7$. In the following section, we will investigate this further.

\section{Dynamic models of repeated, simplified, one-third street Kuhn poker (SKP)}\label{sec_rep}
In SKP, as defined in the previous section, there are three distinct equilibrium solutions. In order to understand which, if any, of these equilibria might be selected in repeated play, we will construct and analyse two models. Each is based on the idea that a player knows their ex-showdown expectation in terms of the three betting frequencies defined above ($b_J$, $c_K$ and $d_Q$) and adjusts the frequency that they control to try to maximise this expectation. This is a reasonable assumption about how human players might approach this game, i.e. bluff more if Players 1 and/or 2 don't call enough, call more if Player 3 bluffs too much. In section~\ref{sec_ode}, we study a model that treats the three betting frequencies as continuous functions of time, $t$, and write down a third order system of nonlinear ordinary differential equations that controls their evolution and has the same three equilibrium solutions as SKP. Although this model has a number of weaknesses, which we discuss, we will see that in some cases its dynamics are very similar to those of the more realistic model that we study in section~\ref{sec_diff}. In particular, none of the three equilibrium states is an asymptotic attractor of the system. The attractors are nested periodic solutions, one set oscillating about the equilibrium that corresponds to Solution 1 and the other about the equilibrium that corresponds to Solution 2, with the selection of the attractor depending on the initial frequencies. In section~\ref{sec_diff}, we describe a difference equation model that is linked to repeated play of SKP, with each player storing information about the most recent plays of the game and using it to estimate the betting frequencies of the other two players. 

Note that in each game the players remain in the same seats for each deal of the cards. In real three player poker games, the participants' roles rotate after every deal, so that in effect they are successively playing three separate games in rotation. Human players make deductions about each others strategies based on their play in similar, but different, situations. We will not consider this possibility here.

\subsection{An ordinary differential equation model of SKP}\label{sec_ode}
The signs of the coefficients of $c_K$, $d_Q$ and $b_J$ in expressions (\ref{eqn_SKP_E1}) to (\ref{eqn_SKP_E3}) respectively, indicate the direction in which each player should change their betting frequency in order to increase their expectation. If we now treat the three frequencies as continuous functions of time, $t$, this suggests that a rational model for how they vary is given by
\begin{equation}
\dot{b}_J = g(b_J) f_3\left(\frac{P-5}{P+1}  - d_Q - c_K + d_Q c_K\right),\label{eqn_ode01}
\end{equation}
\begin{equation}
\dot{c}_K = g(c_K) f_1\left(b_J - \frac{2}{P}\right),\label{eqn_ode02}
\end{equation}
\begin{equation}
\dot{d}_Q = g(d_Q) f_2\left(\frac{c_K-2}{P+1}  +b_J(1-c_K)\right),\label{eqn_ode03}
\end{equation}
where a dot denotes $d/dt$, $g$ is a smooth, non-negative function that vanishes at zero and one, thereby ensuring that each frequency lies in $[0,1]$, and $f_i$  for $i = 1$, $2$, $3$ are non-decreasing functions of their single arguments. The simplest form of this model, which we shall study here, has $g(x) = x(1-x)$ and $f_i(x) = k_i x$, with $k_i$ strictly positive constants. This leads to the system
\begin{equation}
\dot{b}_J =k_3 b_J(1-b_J)\left(\frac{P-5}{P+1}  - d_Q - c_K + d_Q c_K\right),\label{eqn_ode1}
\end{equation}
\begin{equation}
\dot{c}_K = k_1c_K(1-c_K) \left(b_J - \frac{2}{P}\right),\label{eqn_ode2}
\end{equation}
\begin{equation}
\dot{d}_Q = k_2d_Q(1-d_Q)\left(\frac{c_K-2}{P+1}  +b_J(1-c_K)\right).\label{eqn_ode3}
\end{equation}
We can see immediately that the three equilibrium solutions of SKP are also equilibrium solutions of this system, and also that the six planes $b_J=0$, $b_J=1$, $c_K=0$, $c_K=1$, $d_Q=0$ and $d_Q=1$ are invariant. We are only interested in the dynamics of the system for each of $b_J$, $c_K$ and $d_Q$ in $[0,1]$, where the solutions remain. 

It is worth emphasising here that this is not an evolutionary game; there is no population whose composition varies with $t$. The three frequencies are controlled by three individual players, who adjust them in response to the frequencies of the other two players, of which they are assumed to have perfect knowledge. 

A local analysis close to the three equilibrium points, $S_1 = (2/(P+1),0,(P-5)/(P+1))$, $S_2 = (2/P,(P-5)/(P+1),0)$ and, for $P>P^*$, $S_3 = (2/P, 2/(P+2),(P^2-5P-12)/P(P+1))$, which correspond to (\ref{eqn_rsol1}) to (\ref{eqn_rsol3}), shows that
\begin{itemize}
\item $S_1$ has a one-dimensional stable manifold and a two dimensional centre manifold. This centre manifold is the plane $c_K = 0$, and here the dynamics are those of a nonlinear centre, i.e. a series of nested limit cycles (periodic solutions).
\item $S_2$ has a one-dimensional manifold that is stable for $P>P^*$ and unstable for $P<P^*$, and a two dimensional centre manifold. This centre manifold is the plane $d_Q = 0$, and here the dynamics are those of a nonlinear centre, i.e. a series of nested limit cycles.
\item $S_3$ has a one-dimensional unstable manifold and a two-dimensional stable manifold. The dynamics on the stable manifold are oscillatory. The stable manifold separates the two domains of attraction of the planes $c_K = 0$, which contains $S_1$ and $d_Q=0$, which contains $S_2$.
\end{itemize}
Moreover, it is possible to show that this description is also qualitatively correct for the more general system, (\ref{eqn_ode01}) to (\ref{eqn_ode03}). We conclude that, generically, the solution is attracted to one of the limit cycles that surrounds either $S_1$ or $S_2$. When $P<P^*$, the solution is attracted to a limit cycle surrounding $S_1$\footnote{The limit cycles that surround $S_2$ are unstable for $P<P^*$, which reflects the fact that $S_2$ is not an equilibrium solution of SKP.}, but for $P>P^*$  the selection depends on the initial conditions (initial betting frequencies). We will focus on the more interesting case, $P>P^*$.

\begin{figure}
 \begin{center}
 \includegraphics[width=\textwidth]{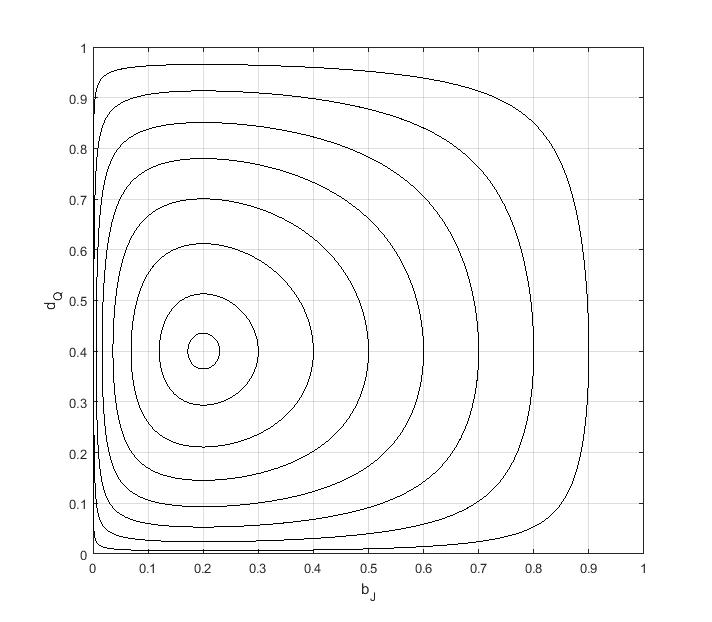} 
 \caption{Nested limit cycles in the plane $c_K=0$.}\label{fig_c0_limitcycles}
 \end{center}
 \end{figure}

\begin{figure}
 \begin{center}
 \includegraphics[width=\textwidth]{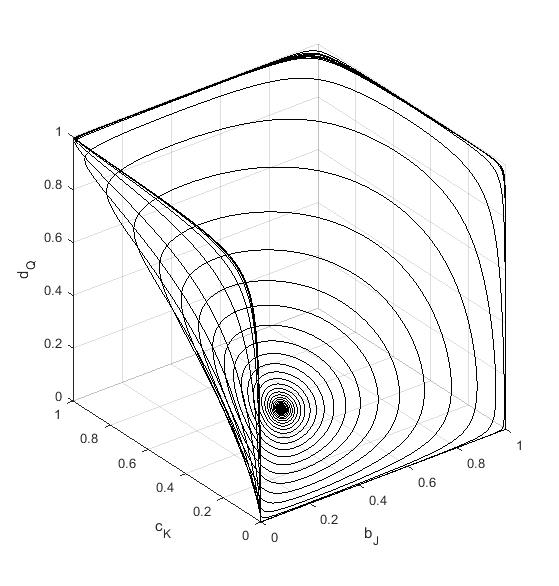} 
 \caption{An integral path in the stable manifold of $S_3$. This manifold separates the basins of attraction of the two attracting planes.}\label{fig_stab_man}
 \end{center}
 \end{figure}

\begin{figure}
 \begin{center}
 \includegraphics[width=\textwidth]{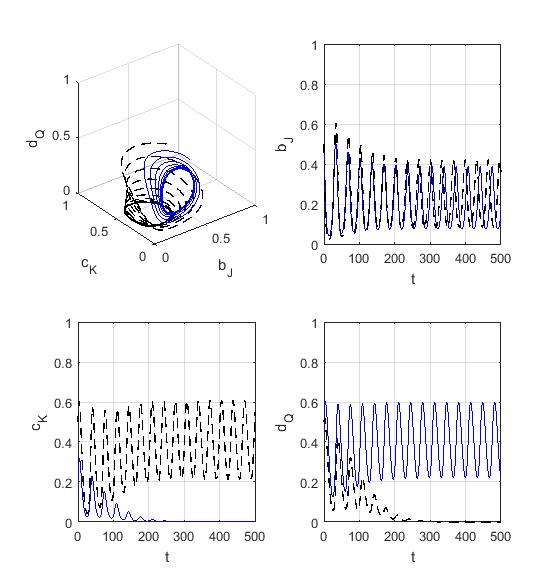} 
 \caption{Two solutions with different final behaviour.}\label{fig_odesols}
 \end{center}
 \end{figure}

\begin{figure}
 \begin{center}
 \includegraphics[width=\textwidth]{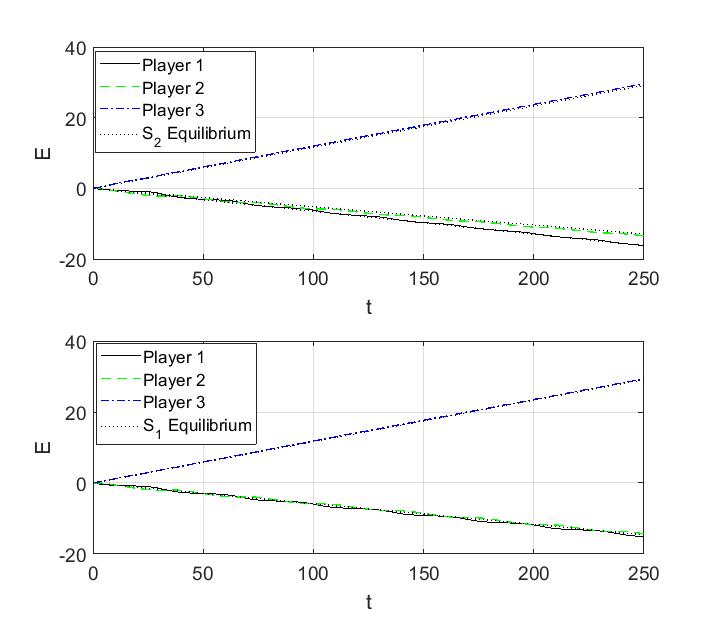} 
 \caption{The ex-showdown expectations of the solutions shown in Figure~\protect{\ref{fig_odesols}}.}\label{fig_odeE}
 \end{center}
 \end{figure}

In the following, all solutions are plotted for the typical case $k_1=k_2=k_3=1$, $P=9$ and calculated numerically using the solver {\it ode45} in MATLAB. Figure~\ref{fig_c0_limitcycles} shows the limit cycles in the plane $c_K=0$. On each integral path, the solution evolves anticlockwise. As Player 3 bluffs more frequently, Player 2 calls more frequently, which causes Player 3 to bluff less frequently, then Player 2 to call less frequently, and so on, {\it ad infinitum}. The solutions in the invariant plane $d_Q=0$ are very similar and represent the same cycle of rise and fall in bluffing and calling frequencies, but for Players 1 and 3. By integrating equations (\ref{eqn_ode1}) to (\ref{eqn_ode3}) backwards in time from an initial point in phase space close to the unstable equilibrium solution $S_3$, we can numerically calculate an integral path that lies in the the stable manifold of $S_2$, which is shown in Figure~\ref{fig_stab_man}. This gives us a visual indication of the boundary of the basins of attraction of the two stable attractors. An example of each type of solution is shown in Figure~\ref{fig_odesols}. 
Finally, Figure~\ref{fig_odeE} shows the ex-showdown expectations of each player for the solutions shown in Figure~\ref{fig_odesols}, defined by
\begin{equation}
E_1(t) = \frac{1}{24}\int_0^t \left\{c_K(s)(Pb_J(s)-2)-(P-2)b_J(s)\right\}\, ds,
\end{equation}
\begin{equation}
E_2(t) = \frac{1}{24}\int_0^t \left\{d_Q(s)\left\{c_K(s)-2+b_J(s)(1-c_K(s))(P+1)\right\} + b_J(s)(c_K(s)+3)-2\right\}\, ds,
\end{equation}
\[
E_3(t) = \frac{1}{24}\int_0^t \left\{b_J(s)\left[P-5 -(P+1)\left\{d_Q(s)+(1-d_Q(s))c_K(s)\right\}\right] \right.
\]
\begin{equation} \left. + 2 c_K(s) +(2-c_K(s))d_Q(s)+2\right\}\, ds.
\end{equation}
We can see that, for the case shown in Figure~\ref{fig_odeE}, these expectations oscillate about those of the $S_1$ and $S_2$ equilibrium solutions given by (\ref{eqn_re1}) and (\ref{eqn_re2}).

\subsection{A difference equation model with estimators and real game play}\label{sec_diff}
In this section, we study a more realistic model of repeated play of SKP. In the $i$th of $N$ rounds of play, Players 1, 2 and 3 have frequencies $c_i$, $d_i$ and $b_i$ respectively and stacks $S_{1i}$, $S_{2i}$ and $S_{3i}$, with $S_{p1} = 0$ for $p = 1$, $2$, $3$. On each round, each player contributes $P/3$ units to the pot, and three cards are dealt, but, for computational efficiency, we exclude combinations where Player 3 is dealt K or Q, since she must check and no decisions are made. One of the twelve possible combinations in which Player 3 is dealt A or J is chosen at random on each round. Play then proceeds as described in section~\ref{sec_simKuhn}, with betting frequencies $b_i$, $c_i$ and $d_i$, and stacks updated appropriately after each round of play is complete, when the betting frequencies are also adjusted using a difference equation analogue of (\ref{eqn_ode1}) to (\ref{eqn_ode3}), namely
\begin{equation}
b_{i+1} = \max \left\{ 0, \min\left\{1, b_i+k_3 \left(\frac{P-5}{P+1}  - \bar{d}_i - \bar{c}_{3i} + \bar{d}_i \bar{c}_{3i}\right)\right\}\right\},\label{eqn_diff1}
\end{equation}
\begin{equation}
c_{i+1} = \max \left\{ 0, \min\left\{1, c_i+k_1 \left(\bar{b}_{1i} - \frac{2}{P}\right)\right\}\right\},\label{eqn_diff2}
\end{equation}
\begin{equation}
d_{i+1} =\max \left\{ 0, \min\left\{1, d_i+ k_2\left(\frac{\bar{c}_{2i}-2}{P+1}  +\bar{b}_{2i}(1-\bar{c}_{2i})\right)\right\}\right\}.\label{eqn_diff3}
\end{equation}
Here, the barred variables on the right hand side are estimators of the opponents' betting frequencies. For example, $\bar{b}_{2i}$ is Player 2's estimate of Player 3's bluffing frequency after $i$ rounds of play. Player $p$ uses the previous $L_p$ hands to construct unbiassed estimators of the other players' frequencies based on information available to him. Details of these estimators are given in Appendix~\ref{sec_estimators}. 

The strategy of each player is characterised by the adjustment rate parameter, $k_p$, which determines how rapidly they adjust their betting frequency in response to their estimates of the opponents' strategies, and $L_p$, the number of recent rounds of play that they use to estimate these strategies. Larger values of $L_p$ lead to more accurate estimates, but a longer delay in the estimates. In this paper, we will focus on the dynamics when each player uses identical parameters, with $L_p = L$ and $k_p = k$ for $p = 1$, $2$, $3$, which allows us to illustrate the dynamic complexity of the game, along with the difficulty of predicting this based on either equilibrium considerations or even the differential equation model studied in the previous section. 

\begin{figure}
 \begin{center}
 \includegraphics[width=\textwidth]{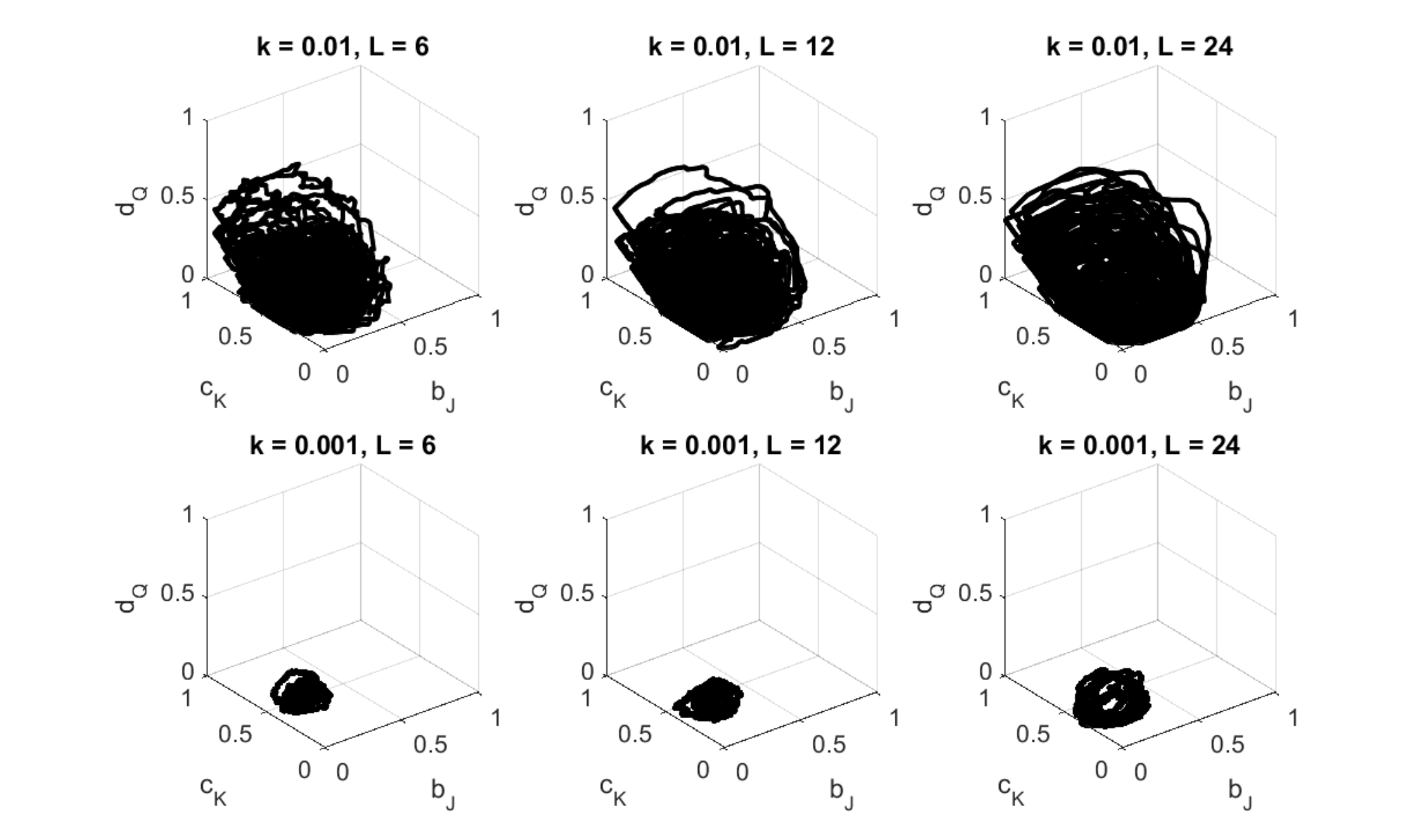} 
 \caption{The solution for $P = 9$, $k = 0.01$ and $0.001$ and $L = 6$, $12$ and $24$.}\label{fig_bcd_P9_1}
 \end{center}
 \end{figure}

\begin{figure}
 \begin{center}
 \includegraphics[width=\textwidth]{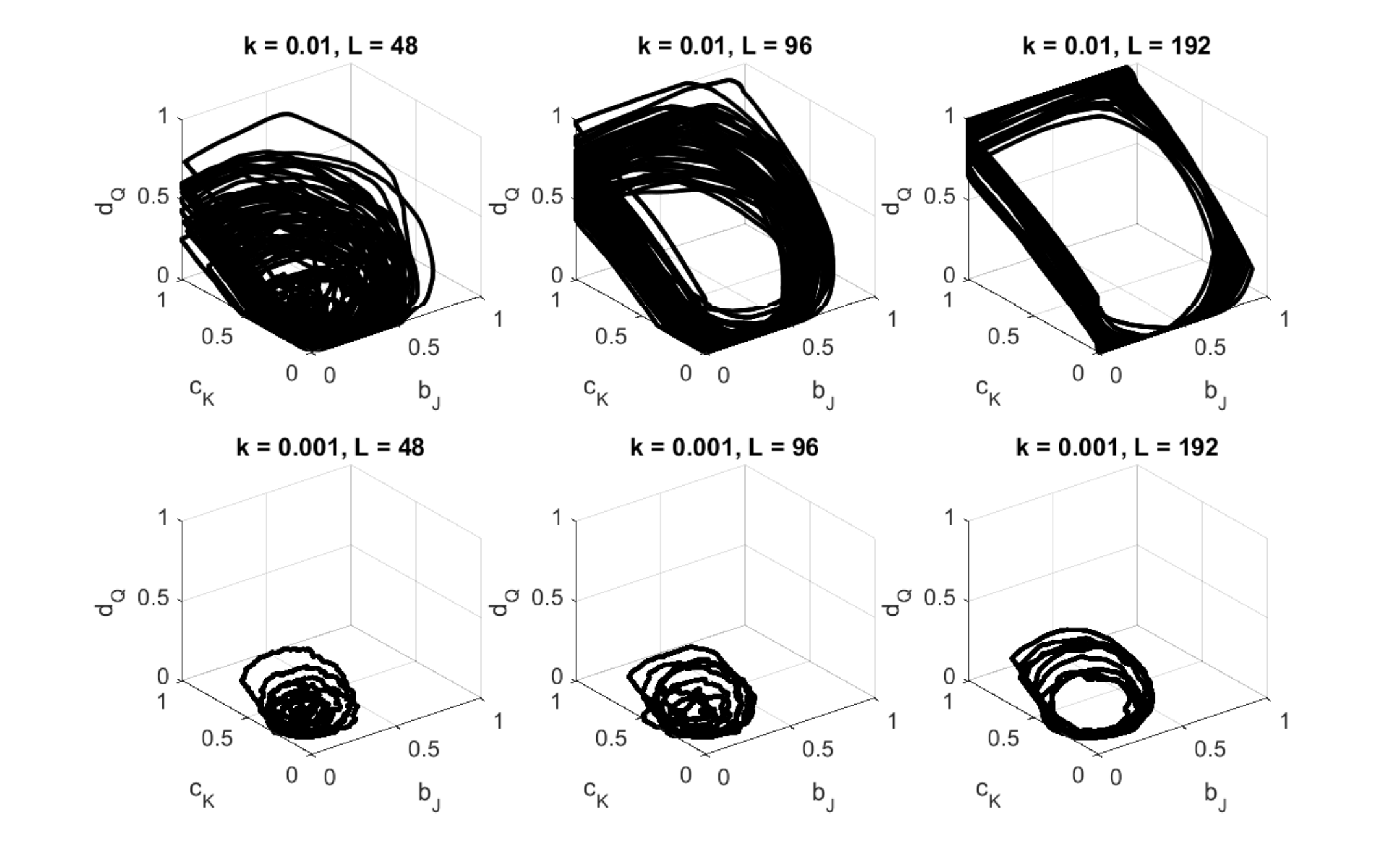} 
 \caption{The solution for $P = 9$, $k = 0.01$ and $0.001$ and $L = 48$, $96$ and $192$.}\label{fig_bcd_P9_2}
 \end{center}
 \end{figure}

\begin{figure}
 \begin{center}
 \includegraphics[width=\textwidth]{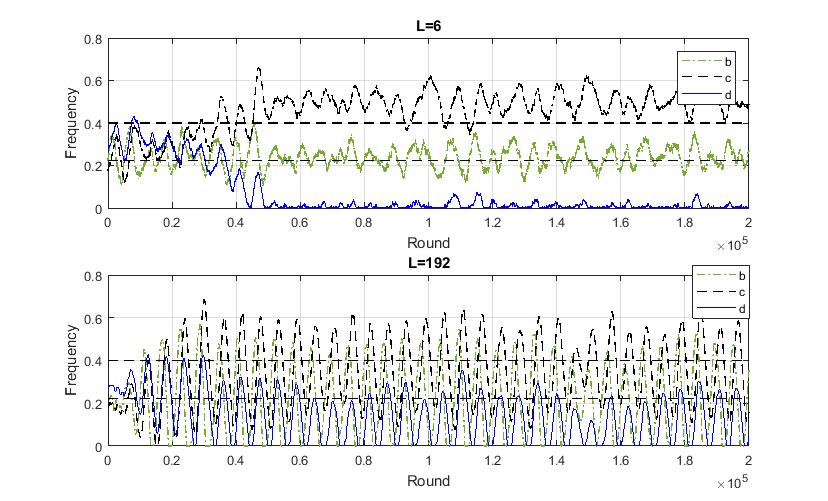} 
 \caption{The solution for $P = 9$, $k=0.001$ and $L = 6$ and $192$. The horizontal dashed lines indicate the solution $S_2$.}\label{fig_P9_L6_L192}
 \end{center}
 \end{figure}

\begin{figure}
 \begin{center}
 \includegraphics[width=\textwidth]{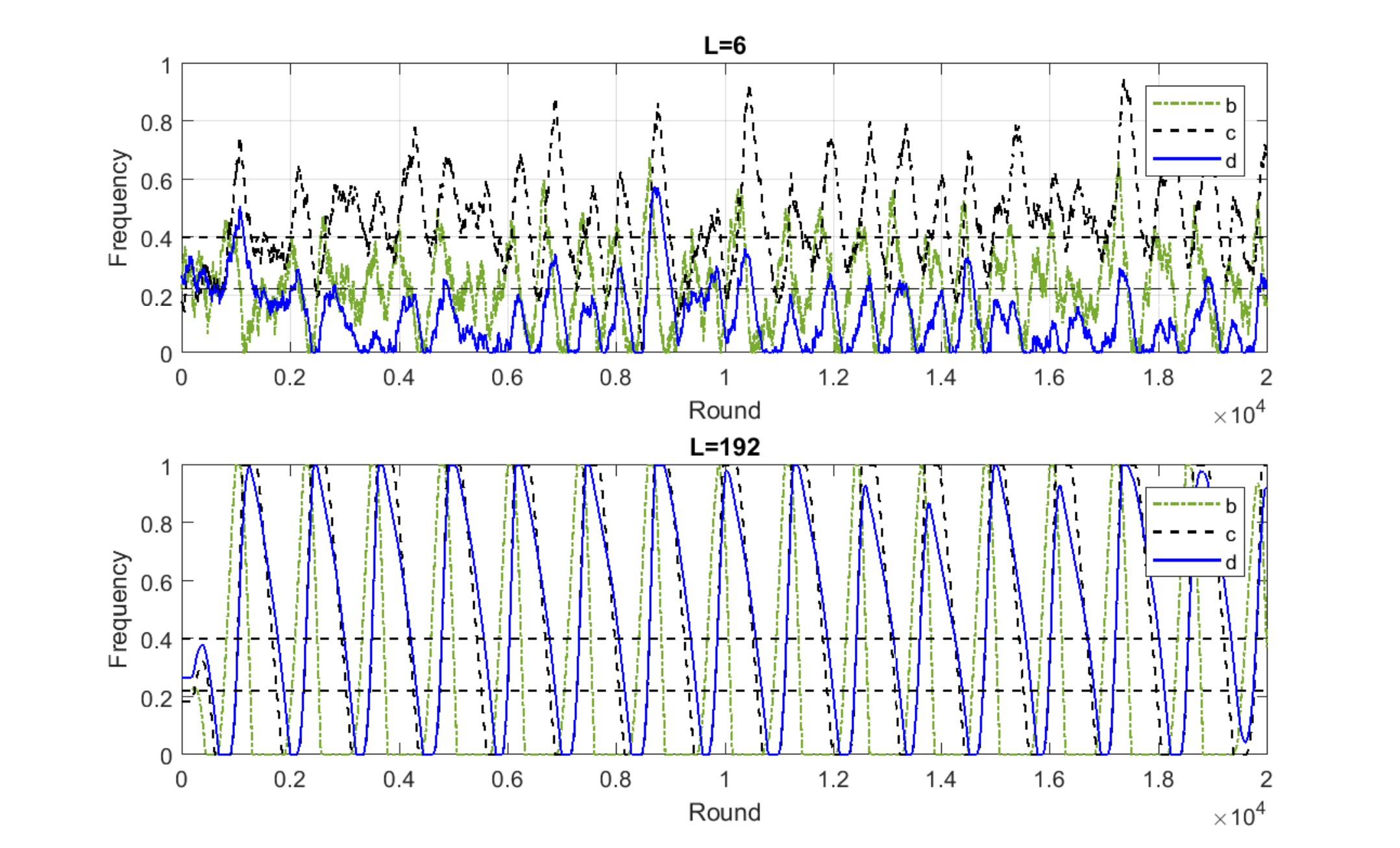} 
 \caption{The solution for $P = 9$, $k=0.01$ and $L = 6$ and $192$. The horizontal dashed lines indicate the solution $S_2$.}\label{fig_P9_L6_L192_2}
 \end{center}
 \end{figure}

\begin{figure}
 \begin{center}
 \includegraphics[width=\textwidth]{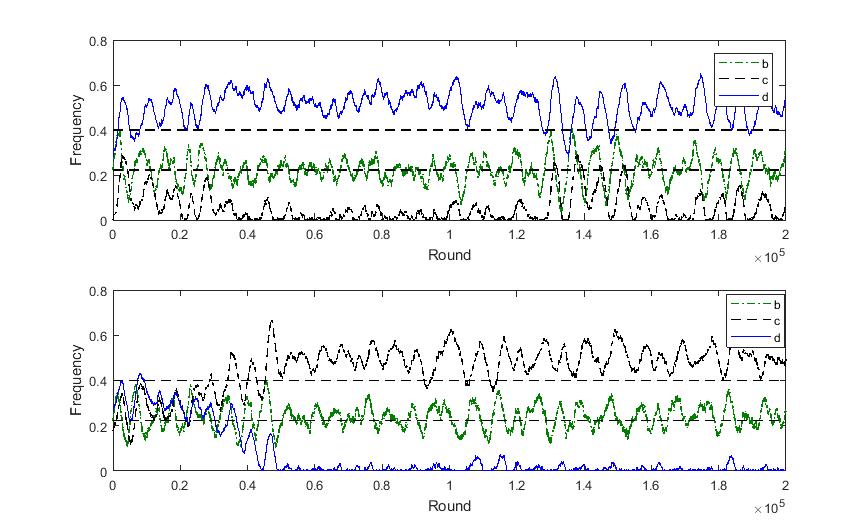} 
 \caption{The solution for $P = 9$, $k=0.001$ and $L = 6$ for two different initial conditions. The horizontal dashed lines indicate the solution $S_2$.}\label{fig_multsol}
 \end{center}
 \end{figure}

Figures~\ref{fig_bcd_P9_1} and~\ref{fig_bcd_P9_2} show how the frequencies vary for a range of values of $k$ and $L$ when $P=9$, with initial conditions close to $S_3$. When the rate of adaptation, $k$, is relatively slow ($k = 0.001$), the solution lies close to $S_2$. Figure~\ref{fig_P9_L6_L192} shows the solution for small and large $L$. Note that the solution becomes smoother as $L$ increases, and lies further from $S_2$. When the rate of adaptation is faster ($k = 0.01$), the solution lies further from $S_2$ and displays rather larger amplitude oscillations, particularly for larger values of $L$. Figure~\ref{fig_P9_L6_L192_2} shows the time series for large and small $L$. The solution when $L = 192$ displays a very regular relaxation oscillation, with Player 3 bluffing more, Players 1 and 2 both calling more in response, Player 3 bluffing less, Players 1 and 2 calling less, and so on. Figure~\ref{fig_multsol} shows that, for $P>P^*$ ($P = 9$ in the Figure), by changing the initial condition, for $k$ and $L$ sufficiently small, the solution can eventually lie close to either $S_1$ or $S_2$. This bistability is a feature that we also saw in solutions of the ordinary differential equation model. 

The ex-showdown expectations of each player in the solutions shown in Figures~\ref{fig_bcd_P9_1} and~\ref{fig_bcd_P9_2} are shown in table~\ref{table_E}. For $k = 0.001$, these are close to those for the equilibrium solution $S_2$, given by (\ref{eqn_re2}), multiplied by a factor of two to account for the fact that only 12 of the 24 possible combinations of cards are dealt in the simulation. For $k = 0.01$, Player 3 does a little better at the expense of Players 1 and 2. 
\begin{table}
\centering
\begin{tabular}{|c| c| c|} 
\hline
~ & $k = 0.001$ & $k = 0.01$\\
\hline 
~ & -0.132 & -0.134\\
$ L = 6$ & -0.098 & -0.101\\
~ & 0.230 & 0.235\\
\hline 
~ & -0.130 & -0.132\\
$ L = 12$ & -0.103 & -0.105\\
~ & 0.234 & 0.237\\
\hline 
~ & -0.129 & -0.133\\
$ L = 24$ & -0.105 & -0.107\\
~ & 0.234 & 0.241\\
\hline 
~ & -0.129 & -0.138\\
$ L = 48$ & -0.105 & -0.110\\
~ & 0.234 & 0.248\\
\hline 
~ & -0.129 & -0.143\\
$ L = 96$ & -0.106 & -0.112\\
~ & 0.235 & 0.255\\
\hline 
~ & -0.130 & -0.139\\
$ L = 192$ & -0.107 & -0.115\\
~ & 0.237 & 0.254\\
\hline 
\end{tabular}
\caption{Ex-showdown expectations when $P = 9$ for Player 1, 2 and 3, corresponding to solutions shown in Figures~\ref{fig_bcd_P9_1} and~\ref{fig_bcd_P9_2}. The corresponding values for Solution 2 are -0.130, -0.104 and 0.233.}
\label{table_E}
\end{table}

\section{Conclusions}\label{sec_conc}
In this paper we have studied a reduced version of three player Kuhn poker - the one-third street game. We found that we could find the complete set of possible equilibrium solutions analytically, for all positive pot sizes $P$. For some values of the pot size ($P=5$, $P > P^* = \frac{1}{2}\left(5+\sqrt{73}\right) \approx 6.77$), there is an opportunity for one player to transfer profit between the other two players at equilibrium. This leads to opportunities for the formation of alliances in repeated play of the game. We also introduced simplified three player, one-third street Kuhn poker (SKP), which has two interesting and useful features. Firstly, there is only one strategic decision available for each player, characterized by a single action frequency. Secondly, we showed that for $P>P^*$, three distinct equilibrium solutions exist, and for $P>7$ there is no obvious way of deciding which solution, if any, is likely to be realized in practice. 

We then moved on to study two dynamic models of repeated play of SKP. In an ordinary differential equation model, we found that the large time solutions are oscillatory, with either Player 1 or Player 2 (but not both) catching bluffs made by Player 3. This model assumes that each player has perfect knowledge of the other players' betting frequencies. In a difference equation model, with each player estimating the other players' betting frequencies based on the information available to them over the previous $L$ deals, we found a variety of possible behaviours, depending on the frequency adjustment rate and memory parameters. For slow enough frequency adjustment, the solutions are noisy, but otherwise similar to the differential equation model, but more rapid frequency adjustment leads to large oscillations in all players' betting frequencies.

\renewcommand{\abstractname}{Acknowledgements}
\begin{abstract}
I would like to acknowledge the contribution of my undergraduate project student at the University of Nottingham in 2015, Richard Farbridge, who discovered two of the three equilibrium solutions in the simplified game.
\end{abstract}

\bibliography{Kuhn}{}
\bibliographystyle{plain}

\begin{appendix}
\section{Kuhn poker: Unique equilibrium solution for $P\geq 2$\label{AKQJ_ap}}
There are $^4\!P_3 = 24$ different ways to deal three from four cards to three players. By considering each of these and the probability and payoffs associated with each possible sequence of actions shown in the game tree in Figure~\ref{fig_tree}, we find that the ex-showdown expectations for each player are given by
\begin{equation}
 24 E_1 = c_Q \left[ -2 + b_J \left\{P-\left(P+2\right)o_K\right\}\right]
+ c_K \left\{ -2 + P \left(b_J+b_Q\right) \right\}
+\left(2-P+o_K\right)\left(b_J+b_Q\right),
\end{equation}
\[
24 E_2 = d_Q \left\{c_K-2+b_J\left(1-c_K\right) \left(P+1\right)\right\} +
d_K\left\{ c_Q - 2 + b_J \left(1-c_Q\right)\left(P+1\right)+b_Q \left(P+1\right)\right\}
\]
\begin{equation}
+ o_K \left[ - c_Q - b_Q + b_J \left\{ \left(P+2 \right) c_Q-1 \right\} \right]
+ b_J \left(2-P+c_Q+c_K\right) + b_Q \left(2-P+c_K\right),
\end{equation}
\[
24 E_3 = b_J \left[ 2P-4 - \left(P+1\right) \left\{ c_Q \left(1-d_K\right) + d_K + c_K \left(1-d_Q \right) + d_Q \right\} \right]
\]
\begin{equation}
+ b_Q \left\{ 2P-4 -\left(P+1\right) \left(c_K+d_K\right) \right\}
+ 2 \left(c_Q+c_K\right) + (2-c_K)d_Q+(2-c_Q)d_K+ c_Q o_K.
\end{equation}
 By noting that Player 3 chooses $b_J$ and $b_Q$ to maximise $E_3$, Player 1 chooses $c_Q$ and $c_K$ to maximise $E_1$ and Player 2 chooses $d_Q$, $d_K$ and $o_K$ to maximise $E_2$, we can find the seven constraints that must be satisfied by equilibrium solutions, which correspond to the seven strategy parameters. In each case there are three possibilities, either an indifference holds (each equality labelled (c)) or the corresponding parameter is chosen to be one or zero to maximise expectation (each inequality labelled (a) or (b)).

 \begin{enumerate}
 \item \begin{enumerate}
 \item $ c_Q(1-d_K)+d_K+c_K(1-d_Q)+d_Q < \frac{2P-4}{P+1} ~~ \mbox{\&} ~~ b_J=1$,
 \item $ c_Q(1-d_K)+d_K+c_K(1-d_Q)+d_Q > \frac{2P-4}{P+1} ~~ \mbox{\&} ~~ b_J=0$,
 \item $ c_Q(1-d_K)+d_K+c_K(1-d_Q)+d_Q = \frac{2P-4}{P+1}~~ \mbox{\&} ~~ 0 \leq b_J \leq 1$,
 \end{enumerate}
 \item \begin{enumerate}
 \item $ c_K + d_K < \frac{2P-4}{P+1} ~~ \mbox{\&} ~~ b_Q=1$,
 \item $ c_K + d_K > \frac{2P-4}{P+1} ~~ \mbox{\&} ~~ b_Q=0$,
 \item $ c_K + d_K = \frac{2P-4}{P+1}~~ \mbox{\&} ~~ 0 \leq b_Q \leq 1$,
 \end{enumerate}
 \item \begin{enumerate}
 \item $ -2 + b_J \left\{ P - \left(P+2 \right)o_K\right\} > 0 ~~ \mbox{\&} ~~ c_Q = 1$,
 \item $ -2 + b_J \left\{ P - \left(P+2 \right)o_K\right\} < 0 ~~ \mbox{\&} ~~ c_Q = 0$,
 \item $ -2 + b_J \left\{ P - \left(P+2 \right)o_K\right\} = 0 ~~ \mbox{\&} ~~ 0\leq c_Q \leq 1$,
 \end{enumerate}
 \item \begin{enumerate}
 \item $ b_J+b_Q > \frac{2}{P} ~~ \mbox{\&} ~~ c_K=1$,
 \item $ b_J+b_Q < \frac{2}{P} ~~ \mbox{\&} ~~ c_K=0$,
 \item $ b_J+b_Q = \frac{2}{P} ~~ \mbox{\&} ~~ 0\leq c_K \leq 1$,
 \end{enumerate}
 \item \begin{enumerate}
 \item $c_K-2 + b_J\left(1-c_K\right) \left(P+1\right) > 0 ~~ \mbox{\&} ~~ d_Q=1$,
 \item $c_K-2 + b_J\left(1-c_K\right) \left(P+1\right) < 0 ~~ \mbox{\&} ~~ d_Q=0$,
 \item $c_K-2 + b_J\left(1-c_K\right) \left(P+1\right) = 0 ~~ \mbox{\&} ~~ 0\leq d_Q \leq 1$,
 \end{enumerate}
 \item \begin{enumerate}
 \item $c_Q-2 + b_J\left(1-c_Q\right) \left(P+1\right) + b_Q\left(P+1\right) > 0 ~~ \mbox{\&} ~~ d_K=1$,
 \item $c_Q-2 + b_J\left(1-c_Q\right) \left(P+1\right)+ b_Q\left(P+1\right)  < 0 ~~ \mbox{\&} ~~ d_K=0$,
 \item $c_Q-2 + b_J\left(1-c_Q\right) \left(P+1\right)+ b_Q\left(P+1\right)  = 0 ~~ \mbox{\&} ~~ 0\leq d_K \leq 1$,
 \end{enumerate}
 \item \begin{enumerate}
 \item $-c_Q - b_Q - b_J + b_J c_Q \left(P+2\right) > 0 ~~ \mbox{\&} ~~ o_K=1$,
 \item $-c_Q - b_Q - b_J + b_J c_Q \left(P+2\right) < 0 ~~ \mbox{\&} ~~ o_K=0$,
 \item $-c_Q - b_Q - b_J + b_J c_Q \left(P+2\right) = 0 ~~ \mbox{\&} ~~ 0\leq o_K \leq 1$.
 \end{enumerate}
 \end{enumerate}

 For each value of $P\geq 2$, our task is to find all sets of strategy parameters, $(b_J, b_Q, c_Q, c_K, d_Q, d_K, o_K)$, such that at least one of (a), (b) and (c) holds for each of these seven constraints. These are the equilibrium solutions.

 We can eliminate 1.(a), 1.(b), 2.(b) and 7.(a) immediately.
 \begin{itemize}
 \item 1.(a) $\implies b_J = 1 \implies 4.(a) \implies c_K = 1 \implies 5.(b) \implies d_Q = 0$. Then 1.(a) $\implies c_K+d_K < \frac{2P-4}{P+1} - c_Q(1-d_K) \leq \frac{2P-4}{P+1} \implies$ 2.(a) $ \implies b_Q = 1 \implies c_K-2 + b_J(1-c_Q) (P+1) + b_Q(P+1) = P+(1-c_Q)(P+1)>0 \implies$ 6.(a) $\implies d_K=1 \implies c_K+d_K > \frac{2P-4}{P+1}$, a contradiction.
 \item 1.(b) $\implies b_J=0 \implies \big(($3.(b) $\implies c_Q = 0)~\&~ ($5.(b) $\implies d_Q = 0)\big)$. Then 1.(b) $\implies c_K+d_K > \frac{2P-4}{P+1} \implies $ 2.(b) $\implies b_Q = 0 \implies \big( (4.(b) \implies c_K = 0) ~\&~ (6.(b) \implies d_K=0)\big) \implies c_Q(1-d_K)+d_K+c_K(1-d_Q)+d_Q = 0$, which contradicts 1.(b).
     \item Since we now know that 1.(c) must hold, $c_K+d_K - \frac{2P-4}{P+1} = -c_Q(1-d_K)-d_Q(1-c_K) \leq 0 \implies$ 2.(b) cannot hold.
 \item 7.(a) $\implies o_K = 1 \implies$ 3.(b) $\implies c_Q = 0$, which contradicts 7.(a).
 \end{itemize}

 The remaining constraints can now be written as
 \begin{enumerate}
 \item \begin{enumerate}\addtocounter{enumii}{2}
\item $c_K+d_K - \frac{2P-4}{P+1} = -c_Q(1-d_K)-d_Q(1-c_K) ~~ \mbox{\&} ~~ 0\leq b_J \leq 1$
\end{enumerate}
 \item \begin{enumerate}
 \item $ -c_Q(1-d_K)-d_Q(1-c_K) < 0 ~~ \mbox{\&} ~~ b_Q=1$,\addtocounter{enumii}{1}
 \item $ -c_Q(1-d_K)-d_Q(1-c_K) = 0 ~~ \mbox{\&} ~~ 0\leq b_Q \leq 1$,
 \end{enumerate}
 \item \begin{enumerate}
 \item $ -2 + b_J \left\{ P - \left(P+2 \right)o_K\right\} > 0 ~~ \mbox{\&} ~~ c_Q = 1$,
 \item $ -2 + b_J \left\{ P - \left(P+2 \right)o_K\right\} < 0 ~~ \mbox{\&} ~~ c_Q = 0$,
 \item $ -2 + b_J \left\{ P - \left(P+2 \right)o_K\right\} = 0 ~~ \mbox{\&} ~~ 0\leq c_Q \leq 1$,
 \end{enumerate}
 \item \begin{enumerate}
 \item $ b_J+b_Q > \frac{2}{P} ~~ \mbox{\&} ~~ c_K=1$,
 \item $ b_J+b_Q < \frac{2}{P} ~~ \mbox{\&} ~~ c_K=0$,
 \item $ b_J+b_Q = \frac{2}{P} ~~ \mbox{\&} ~~ 0\leq c_K \leq 1$,
 \end{enumerate}
 \item \begin{enumerate}
 \item $c_K-2 + b_J\left(1-c_K\right) \left(P+1\right) > 0 ~~ \mbox{\&} ~~ d_Q=1$,
 \item $c_K-2 + b_J\left(1-c_K\right) \left(P+1\right) < 0 ~~ \mbox{\&} ~~ d_Q=0$,
 \item $c_K-2 + b_J\left(1-c_K\right) \left(P+1\right) = 0 ~~ \mbox{\&} ~~ 0\leq d_Q \leq 1$,
 \end{enumerate}
 \item \begin{enumerate}
 \item $c_Q-2 + b_J\left(1-c_Q\right) \left(P+1\right) + b_Q\left(P+1\right) > 0 ~~ \mbox{\&} ~~ d_K=1$,
 \item $c_Q-2 + b_J\left(1-c_Q\right) \left(P+1\right)+ b_Q\left(P+1\right)  < 0 ~~ \mbox{\&} ~~ d_K=0$,
 \item $c_Q-2 + b_J\left(1-c_Q\right) \left(P+1\right)+ b_Q\left(P+1\right)  = 0 ~~ \mbox{\&} ~~ 0\leq d_K \leq 1$,
 \end{enumerate}
 \item \begin{enumerate}\addtocounter{enumii}{1}
 \item $-c_Q - b_Q - b_J + b_J c_Q \left(P+2\right) < 0 ~~ \mbox{\&} ~~ o_K=0$,
 \item $-c_Q - b_Q - b_J + b_J c_Q \left(P+2\right) = 0 ~~ \mbox{\&} ~~ 0\leq o_K \leq 1$.
 \end{enumerate}
 \end{enumerate}

 We will subdivide the remaining analysis according to the value of $c_Q$.

 \subsection{$c_Q = 0$}
 By looking for equilibrium solutions with $c_Q = 0$, after noting that this implies 7.(b), and hence $o_K=0$, the constraints are greatly simplified to
 \begin{enumerate}
 \item\begin{enumerate}\addtocounter{enumii}{2}
\item $c_K+d_K - \frac{2P-4}{P+1}= -d_Q(1-c_K)$, $0\leq b_J \leq 1$,
\end{enumerate}
 \item \begin{enumerate}
 \item $ -d_Q(1-c_K) < 0 ~~ \mbox{\&} ~~ b_Q=1$,\addtocounter{enumii}{1}
 \item $ -d_Q(1-c_K) = 0 ~~ \mbox{\&} ~~ 0\leq b_Q \leq 1$,
 \end{enumerate}
\item\begin{enumerate}\addtocounter{enumii}{1}
 \item $b_J \leq \frac{2}{P}$,
\end{enumerate}
 \item \begin{enumerate}
 \item $ b_J+b_Q > \frac{2}{P} ~~ \mbox{\&} ~~ c_K=1$,
 \item $ b_J+b_Q < \frac{2}{P} ~~ \mbox{\&} ~~ c_K=0$,
 \item $ b_J+b_Q = \frac{2}{P} ~~ \mbox{\&} ~~ 0\leq c_K \leq 1$,
 \end{enumerate}
 \item \begin{enumerate}
 \item $c_K-2 + b_J\left(1-c_K\right) \left(P+1\right) > 0 ~~ \mbox{\&} ~~ d_Q=1$,
 \item $c_K-2 + b_J\left(1-c_K\right) \left(P+1\right) < 0 ~~ \mbox{\&} ~~ d_Q=0$,
 \item $c_K-2 + b_J\left(1-c_K\right) \left(P+1\right) = 0~~ \mbox{\&} ~~ 0\leq d_Q \leq 1$,
 \end{enumerate}
 \item \begin{enumerate}
 \item $b_J+ b_Q> \frac{2}{P+1} ~~ \mbox{\&} ~~ d_K=1$,
 \item $b_J+ b_Q  < \frac{2}{P+1}  ~~ \mbox{\&} ~~ d_K=0$,
 \item $b_J+ b_Q  = \frac{2}{P+1} ~~ \mbox{\&} ~~ 0\leq d_K \leq 1 $.
 \end{enumerate}
 \end{enumerate}
 If 2.(a) is true then $b_Q=1 \implies$ 4.(a) $\implies c_K=1$, which contradicts 2.(a). We conclude that only 2.(c) can be true, and therefore that $d_Q = 0$ or $c_K=1$. Now note from 5.(b) that $c_K = 1 \implies d_Q = 0$, so the least restrictive conclusion is that $d_Q = 0$, and hence from 1.(c)
 \begin{equation}
 c_K+d_K = \frac{2P-4}{P+1},~0\leq b_J \leq 1.\label{eqn1}
 \end{equation}
 Now note that 4.(a) $\implies$ 6.(a) $ \implies c_K+d_K=2$, which contradicts (\ref{eqn1}), so 4.(a) cannot hold. Similarly 6.(b) $\implies$ 4.(b) $\implies c_K+d_K=0$, which also contradicts (\ref{eqn1}), so 6.(b) cannot hold. Then, since 4.(c) and 6.(c) cannot hold simultaneously, we have either
 \begin{equation}
 b_J+b_Q = \frac{2}{P+1},~0\leq d_K \leq 1,~c_K=0,\label{eqn2}
 \end{equation}
 or
 \begin{equation}
 b_J+b_Q = \frac{2}{P},~0\leq c_K \leq 1,~d_K=1.\label{eqn3}
 \end{equation}
 From (\ref{eqn1}) we can see that (\ref{eqn2}) can hold if $2 \leq P \leq 5$, whilst (\ref{eqn3}) can hold if $P\geq 5$, although in this case there is a constraint on the size of $b_Q$ that arises from 5.(b) (shown in (\ref{eq4a})). This is the equilibrium solution given by (\ref{eq1}) to (\ref{eq5}) with $c_Q=0$.

 \subsection{$0<c_Q < 1$}
In this case, 3.(c) shows that
\begin{equation}
b_J = \frac{2}{P-(P+2)o_K} \geq \frac{2}{P},\label{bJeqn}
\end{equation}
and hence 4.(b) cannot hold. It is now easiest to consider 2.(a) and 2.(c) separately.
\subsubsection{$b_Q = 1$}
Now $b_Q = 1 \implies$2.(a)$ \implies \big(($4.(a) $\implies c_K=1)~\&~($6.(a)$\implies d_k = 1)\big)$. However, $c_K = d_K = 1$ is inconsistent with 1.(c), so there is no equilibrium that satisfies 2.(a).
\subsubsection{$ b_Q  < 1$}
Now $b_Q < 1 \implies$ 2.(c)$ \implies \left( d_K = 1~\mbox{and}~d_Q=0 \right)\implies c_K = \frac{P-5}{P+1} \implies (P> 5 ~\&~$4.(c)$) \implies b_J+b_Q = \frac{2}{P}$. Now (\ref{bJeqn}) $\implies b_J = \frac{2}{P},~b_Q = 0~\mbox{and}~o_K=0$. This means that the inequalities in 5.(b), 6.(a) and 7.(b) must hold. On substituting these values of the parameters into these inequalities, we find that we also need
\begin{equation}
P >  P^* \equiv \frac{1}{2}\left(5 + \sqrt{73}\right) \approx 6.77,~~~c_Q < \frac{2}{P+4}.\label{cQposeqn}
\end{equation}
This is Solution B in (\ref{eq5}).

 \subsection{$c_Q=1$}
 If 2.(a) is true $\implies b_Q = 1 \implies \big( ($4.(a) $\implies c_K = 1)~\&~($6.(a) $\implies d_K = 1)\big)$, which contradicts 2.(a), and hence 2.(c) must hold, i.e. $ -(1-d_K)-d_Q(1-c_K) = 0$. Since 5.(b) shows that $c_K = 1 \implies d_Q = 0$, the least restrictive assumption is that $d_Q=0$ and $d_K=1$, and hence $c_K = \frac{P-5}{P+1}$ and $P \geq 5$.

If $o_K = 0$, 3.(a) contradicts 4.(c), so 7.(c) must hold. Along with 4.(c), this gives $b_J = 1/P$ and $b_Q = 1/P$, in contradiction with 3.(a). We conclude that no equilibrium solution with $c_Q = 1$ is possible.

\section{Simplified Kuhn poker: Two or three equilibrium solutions\label{restAKQJ_ap}}
The ex-showdown expectations for this game are given by
\begin{equation}
24E_1 = c_K(Pb_J-2)-(P-2)b_J,\label{eqn_SKP_E1}
\end{equation}
\begin{equation}
24E_2 = d_Q\left\{c_K-2+b_J(1-c_K)(P+1)\right\} + b_J(c_K+3)-2,\label{eqn_SKP_E2}
\end{equation}
\begin{equation}
24E_3 = b_J\left[P-5 -(P+1)\left\{d_Q+(1-d_Q)c_K\right\}\right] + 2 c_K +(2-c_K)d_Q+2,\label{eqn_SKP_E3}
\end{equation}
and the three constraints that must hold at equilibrium are
 \begin{enumerate}
 \item \begin{enumerate}
 \item $ d_Q+(1-d_Q)c_K < \frac{P-5}{P+1}~~ \mbox{\&} ~~ b_J=1$,
 \item $ d_Q+(1-d_Q)c_K > \frac{P-5}{P+1} ~~ \mbox{\&} ~~ b_J=0$,
 \item $ d_Q+(1-d_Q)c_K = \frac{P-5}{P+1} ~~ \mbox{\&} ~~ 0\leq b_J \leq 1$,
 \end{enumerate}
 \item \begin{enumerate}
 \item $ b_J> \frac{2}{P}~~ \mbox{\&} ~~ c_K=1$,
 \item $ b_J< \frac{2}{P} ~~ \mbox{\&} ~~ c_K=0$,
 \item $ b_J = \frac{2}{P} ~~ \mbox{\&} ~~ 0\leq c_K \leq 1$,
 \end{enumerate}
 \item \begin{enumerate}
 \item $ (1-c_K)\left\{b_J(P+1)-1\right\}-1>0~~ \mbox{\&} ~~ d_Q=1$,
 \item $ (1-c_K)\left\{b_J(P+1)-1\right\}-1<0 ~~ \mbox{\&} ~~ d_Q=0$,
 \item $ (1-c_K)\left\{b_J(P+1)-1\right\}-1=0 ~~ \mbox{\&} ~~ 0\leq d_Q \leq 1$.
 \end{enumerate}
 \end{enumerate}

 Firstly, we note that
 \begin{itemize}
 \item 1.(a) $\implies b_J=1 \implies$ 2.(a) $\implies c_K=1 \implies$ 3.(b) $\implies d_Q=0$, in contradiction with 1.(a),
 \item 1.(b) $\implies b_J=0 \implies$ 2.(b) $\implies c_K=0 \implies$ 3.(b) $\implies d_Q=0$, in contradiction with 1.(b).
 \end{itemize}
 We conclude that only 1.(c) can hold, and hence that equilibrium solutions have
 \begin{equation}
  d_Q+(1-d_Q)c_K = \frac{P-5}{P+1}.\label{eqn_r1}
 \end{equation}
 This immediately shows that $c_K \not = 1$ and $d_Q \not = 1$, eliminating the possibility of either 2.(a) or 3.(a) at equilibrium. This leaves just 2.(b), (c) and 3.(b), (c). Although 2.(b) and 3.(b) cannot hold simultaneously, the other three combinations are all possible, and lead to the equilibrium solutions (\ref{eqn_rsol1}) to (\ref{eqn_rsol3}).

\section{Solution using symbolic algebra}
For the relatively small games that we have studied in this paper, the analytical solution can also be determined using symbolic algebra. 

\subsection{Simplifed Kuhn Poker}\label{sec_Math1}
We used Mathematica to confirm the analysis of section~\ref{restAKQJ_ap}. Noting that \&\& and $||$ are the logical AND and OR operators, the command
\begin{verbatim}
FullSimplify[
Solve[((d + (1 - d) c < (P - 5)/(P + 1) &&  b == 1) || 
	(d + (1 - d) c > (P - 5)/(P + 1) && b == 0) || (d + (1 - d) c == (P - 5)/(P + 1))) 
			&&
	 ((b > 2/P &&  c == 1) || (b < 2/P && c == 0) || (b == 2/P)) 
			&& 
	(((1 - c) (b (P + 1) - 1) - 1 > 0 && d == 1) || 
	((1 - c) (b (P + 1) - 1) - 1 < 0 && d == 0) || ((1 - c) (b (P + 1) - 1) - 1 == 0)) 
			&& 
	(d >= 0) && (d <= 1) && (c >= 0) && (c <= 1) && (b >= 0) && (b <= 1)
			&&
	(P>5) , {b, c, d}]]
\end{verbatim}
asks Mathematica to find solutions of the problem defined in section~\ref{restAKQJ_ap}. The resulting solution is
 \begin{center}
 \includegraphics[width=\textwidth]{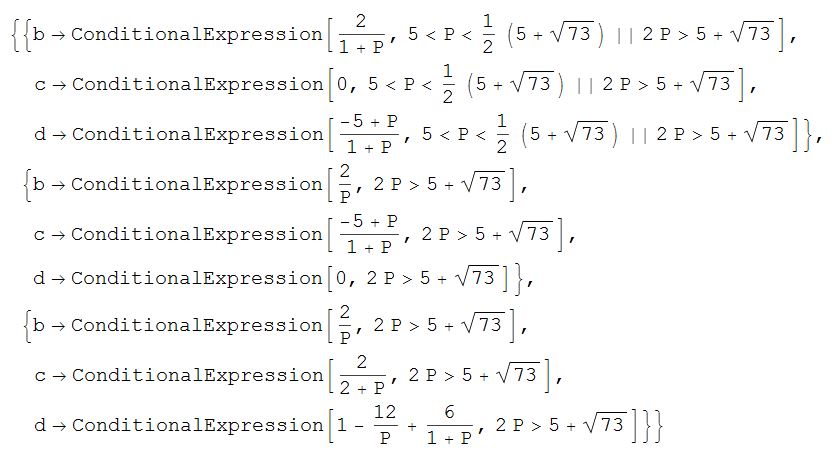}
 \end{center}
which reproduces the solutions (\ref{eqn_rsol1}) to (\ref{eqn_rsol3}) for $P \not = P^*$. In order to find the solution for $P = P^*$, which Mathematica is unable to locate in this general setting, the value of $P$ must be specified and the command run again.

\subsection{Kuhn poker}\label{sec_Math2}
Using this method directly for the problem defined in section~\ref{AKQJ_ap} outruns the 16Gb of RAM available to us. This suggested that we should set up the problem in separate  pieces, using
\begin{verbatim}
S1[1] = cQ (1 - dK) + dK + cK (1 - dQ) + dQ < (2 P - 4)/(P + 1) &&  bJ == 1;
S1[2] = cQ (1 - dK) + dK + cK (1 - dQ) + dQ > (2 P - 4)/(P + 1) &&  bJ == 0;
S1[3] = cQ (1 - dK) + dK + cK (1 - dQ) + dQ == (2 P - 4)/(P + 1);
S2[1] = cK + dK < (2 P - 4)/(P + 1) && bQ == 1;
S2[2] = cK + dK > (2 P - 4)/(P + 1) && bQ == 0;
S2[3] = cK + dK == (2 P - 4)/(P + 1);
S3[1] = -2 + bJ (P - (P + 2) oK) > 0 && cQ == 1;
S3[2] = -2 + bJ (P - (P + 2) oK) < 0 && cQ == 0;
S3[3] = -2 + bJ (P - (P + 2) oK) == 0;
S4[1] = bJ + bQ > 2/P && cK == 1;
S4[2] = bJ + bQ < 2/P && cK == 0;
S4[3] = bJ + bQ == 2/P;
S5[1] = cK - 2 + bJ (1 - cK) (P + 1) > 0 && dQ == 1;
S5[2] = cK - 2 + bJ (1 - cK) (P + 1) < 0 && dQ == 0;
S5[3] = cK - 2 + bJ (1 - cK) (P + 1) == 0;
S6[1] = cQ - 2 + bJ (1 - cQ) (P + 1) + bQ (P + 1) > 0 && dK == 1;
S6[2] = cQ - 2 + bJ (1 - cQ) (P + 1) + bQ (P + 1) < 0 && dK == 0;
S6[3] = cQ - 2 + bJ (1 - cQ) (P + 1) + bQ (P + 1) == 0;
S7[1] = -cQ - bQ - bJ + bJ cQ (P + 2) > 0 && oK == 1; 
S7[2] = -cQ - bQ - bJ + bJ cQ (P + 2) < 0 && oK == 0;
S7[3] = -cQ - bQ - bJ + bJ cQ (P + 2) == 0;
\end{verbatim}
We then ask Mathematica to look for a solution for each of the $3^7 = 2187$ possible combinations of constraints and output every valid solution using
\begin{verbatim}
Do[
sol = FullSimplify[
   Solve[S1[i] && S2[j] && S3[k] && S4[l] && S5[m] && S6[n] &&  S7[o] &&
     (dK >= 0) && (dK <= 1) && (cK >= 0) && (cK <= 1) && (bJ >= 
       0) && (bJ <= 1) && (dQ >= 0) && (dQ <= 1) && (cQ >= 
       0) && (cQ <= 1) && (bQ >= 0) && (bQ <= 1) && (oK >= 
       0) && (oK <= 1) && (P > 0), {bJ, bQ, cK, cQ, dK, dQ, oK}]];
 		If[Length[sol] > 0, {Print[i, j, k, l, m, n, o, sol]}]
 , {i, 1, 3}, {j, 1, 3}, {k, 1, 3}, {l, 1, 3}, {m, 1, 3}, {n, 1, 3}, {o, 1, 3}
 ]
\end{verbatim}
The resulting output is
 \begin{center}
 \includegraphics[width=\textwidth]{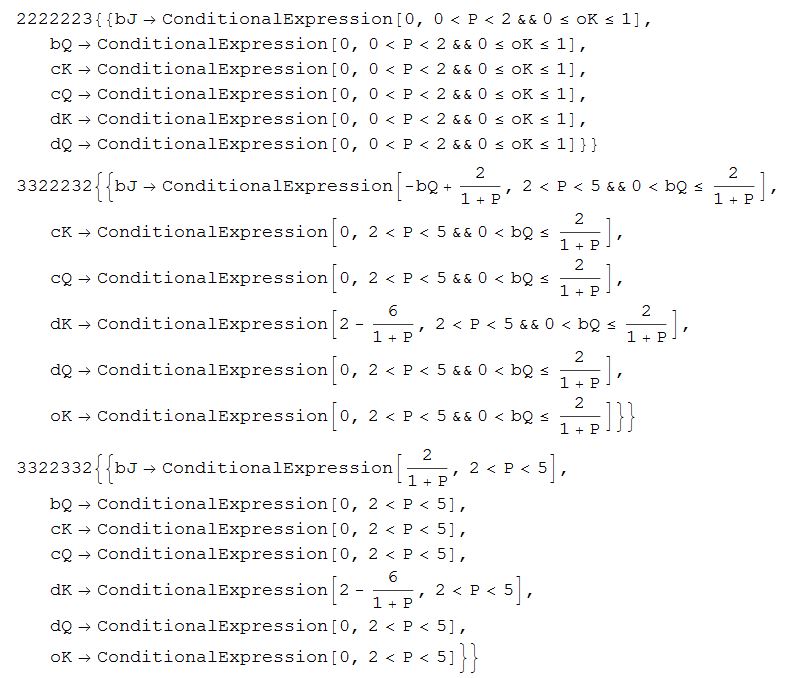}
 \end{center}
 \begin{center}
 \includegraphics[width=\textwidth]{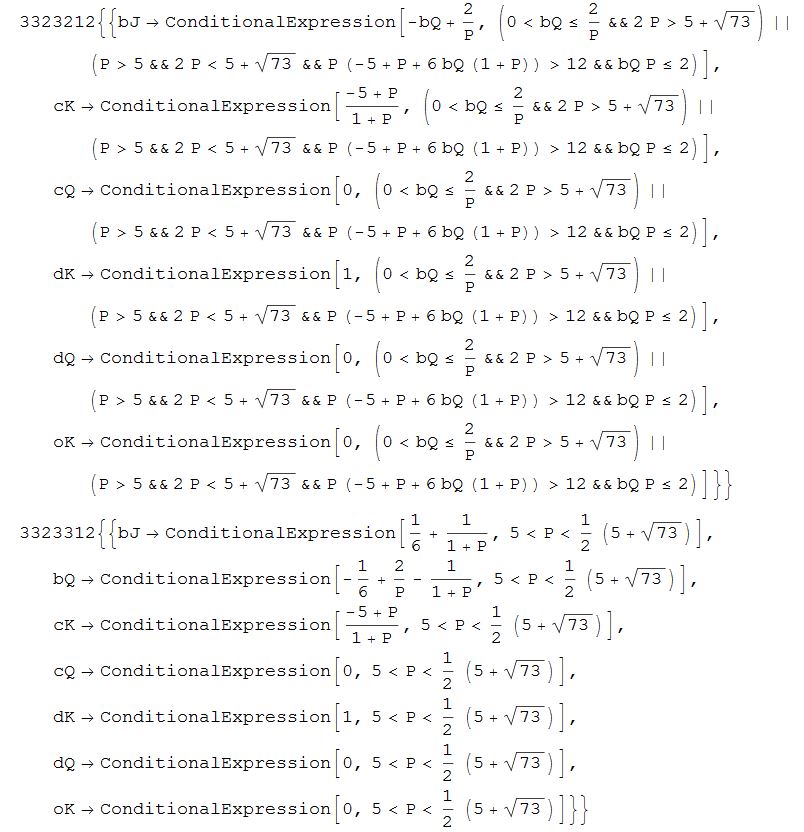}
 \end{center}
 \begin{center}
 \includegraphics[width=\textwidth]{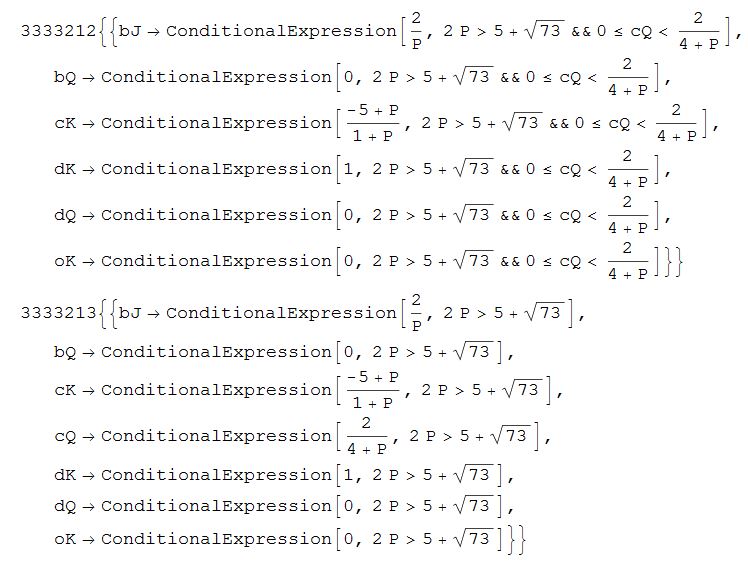}
 \end{center}
This reproduces the solutions (\ref{eq1}) to (\ref{eq5}). Note that, as we found in the previous section, the solutions at the bifurcation points $P=5$ and $P=P^*$ do not appear in this solution and need to be computed separately.

\section{Unbiassed estimators of opponents' betting frequencies}\label{sec_estimators}
In order to describe the estimators used in our simulation, we need to define some notation. Of the twelve possible deals used in the simulation (recall that if Player 3 holds K or Q they check, so no decisions are made and we avoid simulating these deals for reasons of computational efficient) the two cases AJK and AQK (Player 3's card listed first, then those of Players 1 and 2) also do not require any decisions (the play goes bet/fold/call), and hence provide no information to any player, but must be included in the simulation to obtain the correct frequencies. Using the same ordering of cards, we refer to the remaining ten meaningful deals as
\[
D_1 = \mbox{AJQ},~~D_2 = \mbox{AQJ},~~D_3 = \mbox{AKQ},~~D_4 = \mbox{AKJ},~~D_5 = \mbox{JAK},
\]
\begin{equation}
D_6 = \mbox{JAQ},~~D_7 = \mbox{JKA},~~D_8 = \mbox{JQA},~~D_9 = \mbox{JKQ},~~D_{10} = \mbox{JQK}.
\end{equation}
We also define the five possible betting sequences as
\begin{equation}
B_1 = \mbox{k},~~B_2 = \mbox{bff},~~B_3=\mbox{bcf},~~B_4 = \mbox{bfc},~~B_5 = \mbox{bcc},
\end{equation}
with check, bet, call and fold represented by 'k', 'b', 'c' and 'f'.
\begin{sidewaystable}
\centering
  \begin{tabular}{ | c| c| c| c| c| c| c| c| c| c| c| }
    \hline
    ~ & $D_1$ & $D_2$ & $D_3$& $D_4$& $D_5$ & $D_6$ & $D_7$ & $D_8$ & $D_9$ & $D_{10}$\\ 
    ~ & AJQ & AQJ & AKQ& AKJ& JAK & JAQ & JKA & JQA & JKQ & JQK\\ \hline
    $B_1$ & $0$ & $0$& $0$& $0$& $1-b_J$ & $1-b_J$ & $1-b_J$ & $1-b_J$ & $1-b_J$ & $1-b_J$\\ 
    k& ~ & ~ & ~& ~& JAK & JAQ & JKA & JQA & JKQ &JQK\\ \hline
    $B_2$ & $1-d_Q$ & 1 & $(1-c_K)(1-d_Q)$& $1-c_K$& $0$ & $0$ & $0$ & $0$ & $b_J(1-c_K)(1-d_Q)$ & $0$\\ 
    bff & XXX & XXX & XXX& XXX& ~ & ~& ~ & ~ & XXX & ~\\ \hline
    $B_3$ & $0$ & $0$ & $c_K$& $c_K$& $b_J$ & $b_J$ & $0$ & $0$ & $b_Jc_K$ & $0$\\ 
    bcf & ~ & ~ & AKX& AKX& JAX & JAX & ~ & ~& JKQ & ~\\ \hline
    $B_4$ & $d_Q$ & $0$ & $(1-c_K)d_Q$& $0$& $0$ & $0$ & $b_J(1-c_K)$ & $b_J$ & $b_J(1-c_K)$ & $b_J$\\ 
    bfc & AXQ & ~& AXQ& ~& ~ & ~& JXA & JXA & JKQ & JQK\\ \hline
    $B_5$ & $0$ & $0$ & $0$& $0$& $0$ & $0$ & $b_Jc_K$ & $0$ & $0$ & $0$\\ 
    bcc & ~ & ~ & ~& ~& ~ & ~ & JKA & ~ & ~ & ~\\ \hline
  \end{tabular}
\captionof{table}{The information available for each combination of deal and betting sequence.}\label{infotable}
\end{sidewaystable}
For each combination of deal and betting sequence, Table~\ref{infotable} shows the probability of its occurence (public information) and the cards that are known to all players. For example, consider the combination $D_3 = \mbox{AKQ}$ and $B_4 = \mbox{bfc}$. Player 3 always bets holding A, Player 1 folds K with probability $1-c_K$ and Player 2 calls with Q with probability $d_Q$, giving an overall probability $(1-c_K)d_Q$. The publicly known cards are given as AXQ because Players 2 and 3 do not know whether Player 1 folded J or K. The same reasoning gives the entries in the rest of Table~\ref{infotable}.

We now define $N^p_{ij}$ to be the number of occurences of betting sequence $i$ and deal $j$ seen by player $p$ in the last $L_p$ hands (recall that Player $p$ only remembers the last $L_p$ deals). In some cases, a player cannot distinguish between some of these numbers. For example, neither Player 2 nor Player 3 can distinguish between $(B_4, D_1)$ and $(B_4, D_3)$, because they have no way of knowing whether Player 1 held J or K. This means that Players 2 and 3 can only count $N^p_{41}+N^p_{43}$. This is taken into account in the analysis below where necessary (for example, $N^p_{41}+N^p_{43}$ appears in (\ref{est_d}) and (\ref{est_c3}) below).

Players 1 and 2 both estimate Player 3's bluffing frequency using
\begin{equation}
\bar{b}_p = \frac{N^p_{3,5}+N^p_{3,6}+N^p_{4,7}+N^p_{5,7}+N^p_{4,8}+N^p_{4,10}}{N^p_{3,5}+N^p_{3,6}+N^p_{4,7}+N^p_{5,7}+N^p_{4,8}+N^p_{4,10}+\frac{5}{6}\left({N^p_{1,5}+N^p_{1,6}+N^p_{1,7}+N^p_{1,8}+N^p_{1,9}+N^p_{1,10}}\right)}.\label{est_b}
\end{equation}
Only Player 3 needs to estimate the frequency with which Player 2 calls with a Q (see (\ref{eqn_diff2})), and we use
\begin{equation}
\bar{d} = \frac{2\left(N^3_{4,1}+N^3_{4,3}\right) + N^3_{4,9}}{N^3_{4,1}+N^3_{4,3} + N^3_{2,1}+ N^3_{2,2}+ N^3_{2,3}+ N^3_{2,4}+ N^3_{2,9}+ N^3_{4,9}}.\label{est_d}
\end{equation}
Both Player 2 and Player 3 need to estimate the frequency with which Player 1 calls with a K. The asymmetry of the information available to them means that they must use different estimators, and we have taken
\begin{equation}
\bar{c}_2= \frac{2\left(N^2_{3,3}+N^2_{3,4}+N^2_{3,9}+N^2_{5,7}\right)}{N^2_{3,3}+N^2_{3,4}+N^2_{3,9}+N^2_{5,7}+2\left(N^2_{2,2}+N^2_{2,4}\right)+N^2_{4,7}+N^2_{4,8}},\label{est_c2}
\end{equation}
\begin{equation}
\bar{c}_3= \frac{2\left(N^3_{3,3}+N^3_{3,4}\right)+\frac{3}{2}\left(N^3_{3,9}+N^3_{5,7}\right)}{N^3_{3,3}+N^3_{3,4}+N^3_{3,9}+N^3_{5,7}+N^3_{2,1}+N^3_{2,2}+N^3_{2,3}+N^3_{2,4}+N^3_{4,1}+N^3_{4,3}+N^3_{2,9}+N^3_{4,7}+N^3_{4,8}+N^3_{4,9}}.\label{est_c3}
\end{equation}
These are not the only possible estimators, but they are reasonable choices for our simulation. We will derive (\ref{est_d}) here. The other estimators follow in a similar manner.

From the information in Table~\ref{infotable}, we can see that
\[N^3_{4,1}+N^3_{4,3} \propto \left(2-c_K\right)d_Q,\]
\[N^3_{2,1}+N^3_{2,2}+N^3_{2,3}+N^3_{2,4} \propto \left(2-c_K\right)\left(2-d_Q\right),\]
\[N^3_{2,9} \propto b_J \left(1-c_K\right)\left(1-d_Q\right),\]
\[N^3_{4,9} \propto b_J \left(1-c_K\right)d_Q.\]
Hence, as $L_3 \to \infty$ 
\[\frac{N^3_{4,1}+N^3_{4,3}}{N^3_{2,1}+N^3_{2,2}+N^3_{2,3}+N^3_{2,4} } \to \frac{d_Q}{2-d_Q},~~\frac{N^3_{4,9}}{N^3_{2,9}} \to \frac{d_Q}{1-d_Q},\]
and therefore both
\[ \frac{2\left(N^3_{4,1}+N^3_{4,3}\right)}{N^3_{4,1}+N^3_{4,3} + N^3_{2,1}+ N^3_{2,2}+ N^3_{2,3}+ N^3_{2,4}}\]
and
\[ \frac{ N^3_{4,9}}{ N^3_{2,9}+ N^3_{4,9}}\]
are unbiassed estimators of $d_Q$, so we can combine them to give (\ref{est_d}).

\begin{figure}
 \begin{center}
 \includegraphics[width=\textwidth]{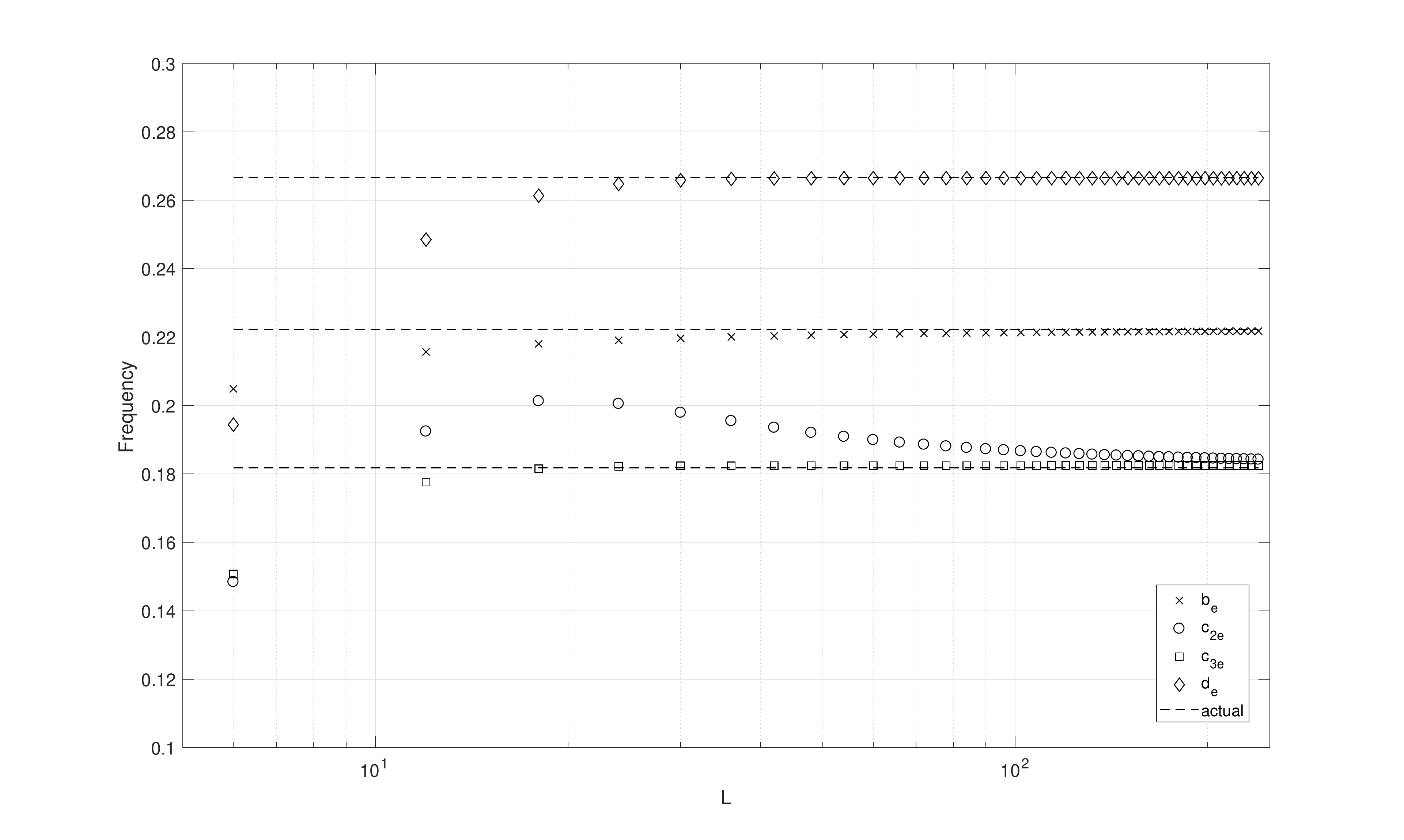} 
 \caption{The mean value of the estimators for $P=9$}\label{fig_meansvL}
 \end{center}
 \end{figure}

\begin{figure}
 \begin{center}
 \includegraphics[width=\textwidth]{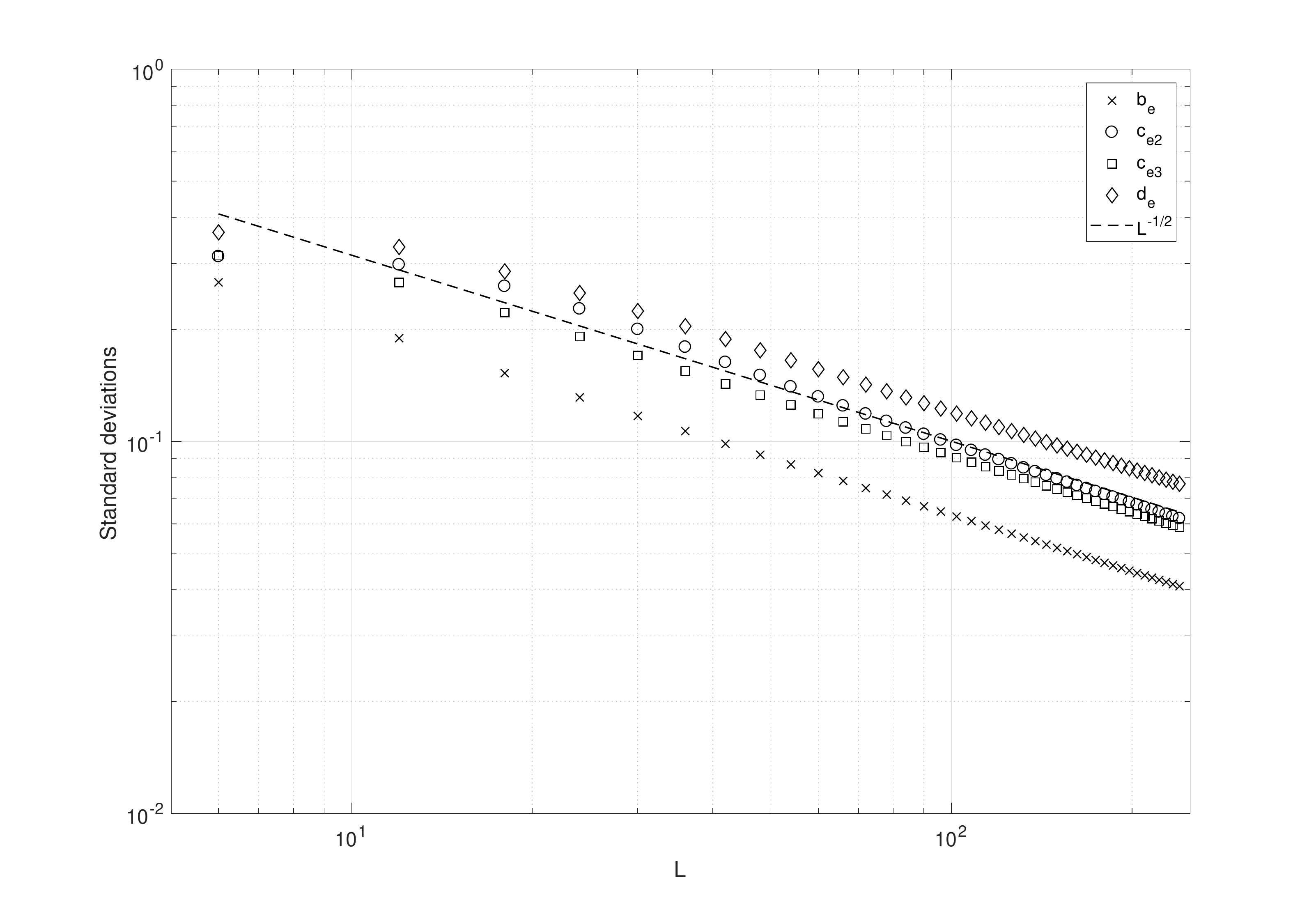} 
 \caption{The standard deviation of the estimators for $P=9$}\label{fig_stdsvL}
 \end{center}
 \end{figure}

In order to confirm that these estimators are unbiassed and to determine how noisy they are as $L$, the number of hands used for the estimate, varies, we simulated $1.1 \times 10^7$ rounds of play with frequencies fixed at Solution 3, given by (\ref{eqn_rsol3}), with $P=9$. We then took the mean and standard deviation of the final $10^7$ values of the estimates. Figure~\ref{fig_meansvL} shows that the estimators converge to the correct frequencies as $L \to \infty$, and indeed that, in the mean, they provide a reasonable estimate of each frequency even when $L$ is as low as six, the smallest value shown. Figure~\ref{fig_stdsvL} shows that the standard deviation of each estimator scales with $L^{-1/2}$ as $L \to \infty$, as we would expect. We can also see that each estimate is very noisy, with standard deviations of all estimates less than 0.1 only for $L>120$, Player 3's estimate of $d$ is the noisiest estimate and the estimate of $b$ the least noisy. This illustrates how hard it is to accurately estimate frequencies without introducing a lag into the estimate. The larger the value of $L$, the more accurate the estimates, at least for constant opponent frequencies, but, for time-varying frequencies, the more out of date the estimate becomes.

\end{appendix}

\end{document}